\documentclass[12pt]{article}
\pdfoutput=1

\usepackage{amsmath}
\usepackage[utf8]{inputenc}
\usepackage[english]{babel}

\usepackage[normalem]{ulem} 
\usepackage{aliascnt}
\usepackage{amssymb}
\usepackage{amsthm}
\usepackage{faktor}
\usepackage{mathtools}
\usepackage{amsfonts}
\usepackage{enumerate}
\usepackage{enumitem}
\usepackage{multicol}
\usepackage{float}
\usepackage{tikz}
\usepackage{tikz-cd} 
\usepackage{comment}
\usepackage{cite}
\usepackage[round,comma,sort,
]{natbib}
\usepackage{hyperref}
\hypersetup{
    colorlinks=true,       
    linkcolor=blue,          
    citecolor=olive,       
    filecolor=magenta,      
    urlcolor=blue          
}

\usepackage[nameinlink]{cleveref}

\DeclareMathOperator{\lk}{Link}
\DeclareMathOperator{\st}{Star}

\newcommand{\cC }{\mathcal C}

\newcommand{\bN}{\mathbb{N}}

\def\coloneqq{\mathrel{\mathop\mathchar"303A}\mkern-1.2mu=}

\begin{document}

\title{On Artin groups admitting retractions \\ to parabolic subgroups}
\author{Bruno A. Cisneros de la Cruz, Mar\'ia Cumplido, Islam Foniqi}

\date{\today}

\maketitle

% ------- Theorem styles -------
\theoremstyle{plain}
\newtheorem{theorem}{Theorem}
\numberwithin{theorem}{section}

\newaliascnt{lemma}{theorem}
\newtheorem{lemma}[lemma]{Lemma}
\aliascntresetthe{lemma}
\providecommand*{\lemmaautorefname}{Lemma}

\newaliascnt{proposition}{theorem}
\newtheorem{proposition}[proposition]{Proposition}
\aliascntresetthe{proposition}
\providecommand*{\propositionautorefname}{Proposition}

\newaliascnt{corollary}{theorem}
\newtheorem{corollary}[corollary]{Corollary}
\aliascntresetthe{corollary}
\providecommand*{\corollaryautorefname}{Corollary}

\newaliascnt{conjecture}{theorem}
\newtheorem{conjecture}[conjecture]{Conjecture}
\aliascntresetthe{conjecture}
\providecommand*{\conjectureautorefname}{Conjecture}

\theoremstyle{remark}

\newaliascnt{claim}{theorem}
\newaliascnt{remark}{theorem}
\newtheorem{claim}[claim]{Claim}
\newtheorem{remark}[remark]{Remark}
\newaliascnt{notation}{theorem}
\newtheorem{notation}[notation]{Notation}
\aliascntresetthe{notation}
\providecommand*{\notationautorefname}{Notation}

\aliascntresetthe{claim}
\providecommand*{\claimautorefname}{Claim}

\aliascntresetthe{remark}
\providecommand*{\remarkautorefname}{Remark}

\newtheorem*{claim*}{Claim}
\theoremstyle{definition}

\newaliascnt{definition}{theorem}
\newtheorem{definition}[definition]{Definition}
\aliascntresetthe{definition}
\providecommand*{\definitionautorefname}{Definition}

\newaliascnt{example}{theorem}
\newtheorem{example}[example]{Example}
\aliascntresetthe{example}
\providecommand*{\exampleautorefname}{Example}

\newaliascnt{question}{theorem}
\newtheorem{question}[question]{Question}
\aliascntresetthe{question}
\providecommand*{\exampleautorefname}{Question}

\begin{abstract}
We generalize the retractions to standard parabolic subgroups for even Artin groups to FC-type Artin groups and other more general families. We prove that these retractions uniquely extend to any parabolic subgroup. We use retractions to generalize the results of Antolín and Foniqi that reduce the problem of intersection of parabolic subgroups to weaker conditions. As a corollary, we characterize coherence for the FC case. 
\end{abstract}

\renewcommand{\thefootnote}{\fnsymbol{footnote}} 
\renewcommand{\thefootnote}{\arabic{footnote}}

%%%%%%%%%%%%%%%%%%%%%%%%%%%%%%%%%%%%%%%%%%%%%%%%%%%%%%%%%%%%%%%%%%%%%%%%%%%%%%%%%%%%%%%%%%%%%%%%%%%%%%%%

% remove notation or backround that is not used in the article

\section{Introduction and notation}
Let~$\Gamma$ be a simple graph  with a finite set of vertices~$S$, and a set of edges~$E$ labelled by elements in~$\{2,3, \dots\} \cup \{ \infty\}$. We define the Artin group associated to~$\Gamma$ as:
$$A_\Gamma \coloneqq \langle S \,|\, \underbrace{sts\cdots}_{m_{s,t}} =\underbrace{tst\cdots}_{m_{s,t}} \text{ whenever } \{s,t\} \in E \text{ and } m_{s,t}\neq \infty  \rangle,$$ 
where~$m_{s,t}$ is the label of the edge joining vertices~$s$ and~$t$. We will refer to~$\Gamma$ as the \emph{Coxeter graph} associated with the Artin group. The definition of Coxeter graphs varies based on their usage: either with no edge if the label is~$\infty$ (no-$\infty$-edge), with no edge if the label is~$2$ (no-$2$-edge), or with all labelled edges (full edge).
In this article, we will use these the three versions of Coxeter graphs, always indicating it beforehand.  When the graph is fixed, we also denote~$A_\Gamma$ by~$A_S$. 

One can also define the \emph{Coxeter group based on~$\Gamma$} which is the group~$A_{\Gamma} / \langle \langle s^2 \mid s \in S \rangle \rangle$. 
%with a presentation 
%$$W_\Gamma \coloneqq \langle S \,|\, s^2 = 1, \, \underbrace{sts\cdots}_{m_{s,t}} =\underbrace{tst\cdots}_{m_{s,t}} \text{ whenever } \{s,t\} \in E \text{ and } m_{s,t}\neq \infty  \rangle.$$ 
An Artin group~$A_\Gamma$ is called of \emph{spherical type} if the associated Coxeter group~$W_\Gamma$ is finite, see for example \citep{Bourbaki, BrieskornSaito}. 

\smallskip

Despite their simple definition and the thorough understanding of Coxeter groups, Artin groups have posed a challenge for the community of group theorists ever since their study began in the 1960s with Jacques Tits and the Bourbaki group \citep{Bourbaki}. 
For example, basic questions such as the word problem or the conjugacy problem still do not have a general answer. However, several families of Artin groups have been found to possess a rich geometric structure: actions on different simplicial complexes \citep*{CharneyDavis,ConciniSalvetti,Paris2014,cumplido2019parabolic} have been discovered, and their study has been very fruitful. One of the main ingredients for the construction of these complexes are the so-called parabolic subgroups. Given~$X\subset S$, the \emph{standard parabolic subgroup}~$A_X$ is the subgroup of~$A_\Gamma$ generated by~$X$. By \cite{van1983homotopy} we know that~$A_X=A_{\Gamma'}$ where~$\Gamma'$ is the subgraph of~$\Gamma$ induced by~$X$. Any conjugate of a standard parabolic subgroup is called a \emph{parabolic subgroup}.

\smallskip

We say that a group morphism~$\rho \colon G\rightarrow H$, where~$H \leq G$, is a \emph{retraction} if~$\rho|_{H}=id_H$. Given an Artin group~$A$ with only even labels, it admits an evident retraction~$\rho: A \rightarrow A_X$ defined by~$\rho(s)=1$ if~$s \not\in X$ and~$\rho(s)=s$ otherwise. 
\cite{AntolinFoniqi} use the already known retractions of even Artin groups to standard parabolic subgroups, and they manage to prove several results in even FC-type Artin groups, including that the intersection of parabolic subgroups is again a parabolic subgroup. 

\smallskip
This article is a first approach to the generalization of the ideas in \citep{AntolinFoniqi} and \citep{antolin2023subgroups} to other classes of Artin groups where we could have odd labels too. Here we will explore how one can define these retractions and their use to obtaining analogous results as in the aforementioned articles.

\smallskip

During the rest of the article, when we say \emph{Artin groups admitting retractions}, we will be referring to Artin groups that admit retractions to standard parabolic subgroups. In \Cref{Sec: FC-type Artin groups admitting retractions}, we give a complete classification of FC-type Artin groups admitting retractions to standard parabolic subgroups. Furthermore, we give them explicitly, and in particular we prove that they send generators to generators or to the identity. In \Cref{sec: Retractions to standard parabolic subgroups uniquely extending to parabolic subgroups} we show that an Artin group admitting retractions to every standard parabolic subgroup also admits retractions to all of its parabolic subgroups. Moreover, for FC-type and other cases, these retractions are uniquely determined for every parabolic subgroup~$P$, i.e., the definition for the retraction does not depend on how we choose to express~$P=\alpha A_X \alpha^{-1}$. In \Cref{sec: Retractions and intersection of parabolic subgroups} we adapt results from \citep{AntolinFoniqi} for even Artin groups, to Artin groups with retractions sending generators to generators or identity and satisfying additional conditions (where FC-type admitting retractions are included). In particular, we reduce the problem of intersection of parabolic subgroups to a weaker condition --- see \autoref{thm: reduction to stars} ---. Finally, in \Cref{sec: Coherence} we also expand the results of \citep{antolin2023subgroups} to characterize coherence in our class of FC-type Artin groups admitting retractions.

\section{FC-type Artin groups admitting retractions}\label{Sec: FC-type Artin groups admitting retractions}

Our first aim is to identify a family of Artin groups admitting retractions, which is different from the family of even Artin groups. In this section, we completely classify the family of FC-type Artin groups admitting retractions 
%to standard parabolic subgroups 
and give explicit retractions for each case. An FC-type Artin group~$A_\Gamma$ is defined by the property that for every subgraph~$\Delta$ of~$\Gamma$ without~$\infty$~edges,~$A_\Delta$ has spherical type --- in other words the associated Coxeter group over~$\Delta$ is finite.
%Spherical-type Artin groups are those that become finite when the relation~$s^2=1$ is added for each generator~$s$ -- see, for example, \citep{Bourbaki, BrieskornSaito}. 
The family we focus on appears naturally, as many well-known results already exist for it, and FC-type Artin groups are also studied in \citep{AntolinFoniqi} through the use of retractions. We show that the family of FC-type Artin groups admitting retractions contains properly the family of even FC-type Artin groups.

\smallskip

In this section, we will solve some equations in dihedral Artin groups. For elements~$x$ and~$y$, we will denote by~${}_s(x,y)$ the alternating product~$xyx\cdots~$ of~$s$ factors, and by~$(x,y)_s$ the alternating product~$\cdots yxy$ of~$s$ factors.

\begin{proposition}\label{prop: eq free group}
Let~$F_2$ be a free group generated by~$a$ and~$b$. For natural numbers~$r,s$ with~$r$ odd, and~$s$ even, the system of equations~$(a,x)_r = (x,a)_r$ and~$(b,x)_s = (x,b)_s$, is incompatible in~$F_2$.
\end{proposition}

To prove the proposition, we will use two well known results about free groups, stated in the following lemma.
\begin{lemma}\label{lemma: some well-known stuff in free group}
Let~$F$ be a free group, and let~$u, v \in F$ be any two non-trivial elements. 
\begin{itemize}
    \item[(i)] If~$u^n = v^n$ for some~$n \in \mathbb{N}$, then~$u = v$.
    \item[(ii)] The centralizer of~$v$ in~$F$ is the cyclic group generated by~$\langle v \rangle$.
\end{itemize}
\end{lemma}

\begin{proof}[Proof of \autoref{prop: eq free group}]
The equation~$(b,x)_s = (x,b)_s$ for an even~$s$, is equivalent to the equation~$(bx)^\frac{s}{2} = (xb)^\frac{s}{2}$, so any solution for it, is also a solution for~$bx = xb$. Now~$x$ belongs to the centralizer of~$b$ so~$x$ is in the cyclic group~$\langle b \rangle$. Assume that~$x = b^k$ for some~$k \in \mathbb{Z}$. Substituting it into~$(a,x)_r = (x,a)_r$ yields an equality of two distinct reduced words in~$F_2$, hence there is no solution to the system.
\end{proof}

\begin{comment}
\begin{proof}[Proof of \autoref{prop: eq free group}]
Suppose that there is a solution~$x$ which is expressed as a word in its reduced normal form written in~$a, b$ and their inverses --- that is, there are no instances in~$x$ of~$dd^{-1}$ for any~$d \in \{a, a^{-1}, b, b^{-1}\}$. If~$x$ is a solution of both equations, it cannot have~$a$ or~$b$ as its first or last letter; because one of the equations~$(a,x)_r = (x,a)_r$ or~$(b,x)_s = (x,b)_s$ would have reduced distinct words on both sides. 
For the same reason,~$x$ cannot start and finish simultaneously with~$a^{-1}$ (or~$b^{-1}$). Therefore,~$x$ needs to have at least 2~letters and must start with~$a^{-1}$ and end with~$b^{-1}$, or conversely start with~$b^{-1}$ and end with~$a^{-1}$. Suppose that we are in the first case (the second case is analogous).

From the second equation, we know that~$(b,x)_{s-1}$ centralizes~$b$, so it must be a power of~$b$. To achieve this, since~$x$ starts with~$a^{-1}$, the second-to-last letter must be~$a$. In the same spirit, if the second letter is~$s$, the third-to-last letter must be~$s^{-1}$, and so on. Therefore,~$x = a^{-1}ww^{-1}a b^{-1}$, which is impossible by assumption.
\end{proof}
\end{comment}

\noindent 
If~$X$ is a three-element subset of the set of standard generators of an Artin group~$A_S$, we call the subgroup~$A_X$ a \emph{triangle} and we will refer to it using the notation~$(m_1,m_2,m_3)$, where the numbers~$m_i$ are the labels of the three corresponding edges formed from elements of~$X$.

\begin{proposition}\label{lem:impossible triangles infinity}
    Let~$m , n$ be positive integers. Artin groups that correspond to the triangles~$(\infty, 2m, 2n+1)$ and~$(\infty, 2m+1, 2n+1)$ do not admit retractions.
\end{proposition}

\begin{proof}
Let us first treat the~$(\infty, 2m, 2n+1)$ case. Let~$A_{a,b,c}$ be an Artin group such that~$m_{a,b} = \infty$,~$m_{a,c} = 2m$, and~$m_{b,c} = 2n+1$. Suppose that there is a retraction~$\rho$ to~$A_{a,b}$ such that~$\rho(c) = x$. The element~$x$ must satisfy the equations~$(a,x)_{2m} = (x,a)_{2m}$~and~$(b,x)_{2n+1} = (x,b)_{2n+1}$ in~$F_2$, and we know by \autoref{prop: eq free group} that there is no solution to this system of equations.

\noindent For the~$(\infty, 2m+1, 2n+1)$ case, notice that any solution of~$(a,x)_{2m+1} = (x,a)_{2n+1}$ is also a solution of~$(a,x)_{4m+2} = (x,a)_{4m+2}$, which together with the equation~$(b,x)_{2n+1} = (x,b)_{2n+1}$ provides an incompatible system of equations in~$F_2$.
\end{proof}

\begin{comment}
\begin{proof}
Let us first treat the~$(\infty, 2k, 2s+1)$ case. Let~$A_{a,b,c}$ be an Artin group such that~$m_{a,b} = \infty$,~$m_{a,c} = 2k$, and~$m_{b,c} = 2s+1$. Suppose that there is a retraction~$\rho$ to~$A_{a,b}$ such that~$\rho(c) = x$. The element~$x$ must satisfy the equations~$(a,x)_{2k} = (x,a)_{2k}$ and~$(b,x)_{2s+1} = (x,b)_{2s+1}$ in~$F_2$, and we know by \autoref{prop: eq free group}that there is no solution to this system of equations.

For the~$(\infty, 2l+1, 2s+1)$ case, notice that any solution of the equation~$(a,x)_{2l+1} = (x,a)_{2l+1}$ is also a solution for~$(a,x)_{4l+2} = (x,a)_{4l+2}$, so again we have an incompatible system of equations.
\end{proof}
\end{comment}

\begin{proposition}\label{prop:characterization of retractions}
Artin groups corresponding to triangles~$(2,3,3)$ and~$(2,3,5)$ do not admit retractions.
\end{proposition}

\begin{proof}
Take first the triangle~$(2,3,3)$ defining an Artin group~$A$ with generators~$\{a,b,c\}$ where~$m_{a,b}=2$,~$m_{b,c}=3$,~$m_{c,a}=3$. Suppose that there is a retraction~$\rho$ from~$A$ onto~$A_{\{a,b\}}$. However, any surjective morphism from~$A$ onto~$A_{\{a,b\}}$ would induce a surjective morphism on their abelianizations,~$\mathbb{Z}$ and~$\mathbb{Z}^2$ respectively; but one cannot have a surjective morphism from~$\mathbb{Z}$ onto~$\mathbb{Z}^2$. 

\noindent 
The same argument works for the triangle~$(2,3,5)$.
\end{proof}

For the next proposition we will use Garside normal forms in the dihedral Artin group with a relation of length 4 and generators~$a$ and~$b$ --- see \citep{BrieskornSaito,DehornoyParis}. In this case, we know that the centre of the group is generated by~$abab$ and that every element~$x$ can be uniquely expressed on its normal form~$(abab)^p x_1 \cdots x_r$, with~$p \in \mathbb{Z}$ and each~$x_i$ equal to~$a$,~$b$,~$ab$,~$ba$,~$aba$ or~$bab$, satisfying the condition that the last letter of~$x_i$ and the first letter of~$x_{i+1}$ must be the same. 

\begin{proposition}\label{prop: characterization of retractions}
The Artin group corresponding to triangle~$(2,3,4)$ does not admit retractions.
\end{proposition}

\begin{proof}
Let~$A$ be an Artin group generated by~$\{a, b, c\}$ with~$m_{a,b}=4$,~$m_{b,c}=3$, and~~$m_{c,a}=2$. Suppose that there is a retraction~$\rho$ from~$A$ onto~$A_{a,b}$ and set~$\rho(c)=x$. Since~$\rho$ is a homomorphism, the equations~$ax=xa$ and~$bxb=xbx$ must be satisfied in~$A_{a,b}$. By \citep[Theorem~5.1]{Paris},~$x$ is a product of powers of~$a$ and positive powers of~$bab$ or~$b^{-1}a^{-1}b^{-1}$.
%--- see \autoref{conj Godelle} for more details. 
When putting this in left normal form, we have two options:~$x = (abab)^s a^t$ or~$x = (abab)^s (bab)^t$, where~$t \geq 0$. Notice that, since~$abab = baba$ in~$A_{\{a,b\}}$, any equality between words demands that the exponent sum of~$b$'s in each word is the same. For any~$x$ as above, the exponent sum of~$b$'s in~$xbx$ is odd and in~$bxb$ is even, leading to a contradiction.
\end{proof}

\begin{theorem}\label{thm: characterization of FC type admitting retractions}
    An FC-type Artin group~$A_S$ admits retractions to parabolic subgroups if and only if its Coxeter graph~$\Gamma$ only contains triangles of the form~$(2,2,k)$,~$(2m,2n,\infty)$, or~$(k,\infty,\infty)$, for~$m, n > 0$, and~$k>1$.

    Moreover, for~$X \subseteq S$ we define the ordinary retraction~$\rho_X: A_S \rightarrow A_X$ as:~$\rho(x)=x$ for all~$x \in X$,~$\rho(v)=a$ if there is an odd edge connecting the generator~$v\in S \setminus X$ to~$a\in X$ and~$\rho(v)=1$ otherwise.
\end{theorem}

\begin{proof}
    If a triangle on~$\Gamma$ contains some~$\infty$, then by \autoref{lem:impossible triangles infinity}, these triangles can only be of the form~$(2m,2n,\infty)$ or~$(k,\infty,\infty)$. Otherwise, the FC condition tells us that the only triangles we can have are~$(2,2,k)$,~$(2,3,3)$,~$(2,3,4)$, and~$(2,3,5)$, which are the only triangles having spherical type. By \autoref{prop: characterization of retractions}, the only triangle admitting retractions to parabolic subgroups from that list is~$(2,2,k)$.

\smallskip
Now we need to prove that any FC-type Artin group whose defining graph contains only triangles given as above, admits retractions. It is obvious that the retraction described in the statement is well-defined for these triangles. Let~$X$ be a proper subset of~$S$. To show that~$\rho$ is a retraction, we need to verify that~$\rho$ preserves the defining relations. Consider an edge~$\{v_1, v_2\}$ of~$\Gamma$ with label~$l$.  We examine the following cases:
\begin{itemize}
    \item[(1)]~$v_1, v_2 \in X$. Here both generators are fixed, so the relation is preserved.
    % one has~$\rho(v_1) = v_1$, and~$\rho(v_2) = v_2$
    \item[(2)]~$v_1 \in S \setminus X$,~$v_2 \in X$ (the other case in analogous). Here we have~$\rho(v_2) = v_2~$. Now we have two options for~$\rho(v_1)$:
    \begin{itemize}
        \item[(2.1)] If~$l$ is odd, then~$\rho(v_1) = v_2$. The relation~${}_l(v_1, v_2) = {}_l(v_2, v_1)$, after applying~$\rho$, becomes~$v_2^l = v_2^l$ in~$A_X$.\item[(2.2)] If~$l$ is even, we have two situations: \\
      ~$\bullet$ If any edge connecting~$v_1$ to~$X$ is even or~$\infty$, then~$\rho(v_1) = 1$. Hence, the relation~${}_l(v_1, v_2) = {}_l(v_2, v_1)$, after applying~$\rho$, becomes~$v_2^{l/2} = v_2^{l/2}$ in~$A_X$. \\
      ~$\bullet$ Suppose that there is an odd edge connecting~$v_1$ to~$X$, say~$\{v_1, x\}$. Then~$l=2$ and the edge~$\{v_2, x\}$ also has label 2. In this case,~$\rho(v_1) = x$, so the relation~$v_1 v_2 = v_2 v_1$, after applying~$\rho$, becomes~$x v_2 = v_2 x$ in~$A_X$, as required.
    \end{itemize}
    \item[(3)]~$v_1, v_2 \in S \setminus X$.
    \begin{itemize}
        \item[(3.1)] Suppose~$l$ is odd; since we cannot have two consecutive odd labels, one has that labels from~$v_1, v_2$ to~$X$ are either equal to~$2$, or to~$\infty$. So~$\rho(v_1) = 1 = \rho(v_2)$. The relation~${}_l(v_1, v_2) = {}_l(v_2, v_1)$ after applying~$\rho$ becomes~$1 = 1$ in~$A_X$.
        \item[(3.2)]If~$l$ is even, we have the following cases: \\
      ~$\bullet$ If~$l > 2$, then both~$v_1$ and~$v_2$ are connected to~$X$ by even labels, so~$\rho(v_1) = 1$ and~$\rho(v_2) = 1$. Hence, the equality of the relation is obviously preserved. \\
      ~$\bullet$ If~$l = 2$ and both~$v_1$ and~$v_2$ are connected to~$X$ by even labels, then again the equality of the relation is preserved as~$\rho(v_1) = 1$ and~$\rho(v_2) = 1$. \\
      ~$\bullet$ Finally, if~$l = 2$ and one of~$v_1$ or~$v_2$, say~$v_1$, is connected to~$x \in X$ by an odd label, then~$\rho(v_1) = x$ and~$\rho(v_2) = 1$. Therefore, the relation~$v_1 v_2 = v_2 v_1$, after applying~$\rho$, becomes~$x = x$ in~$A_X$.
    \end{itemize}
\end{itemize}
\end{proof}

\begin{remark}
Notice that in FC-type Artin groups admitting retractions, odd-labeled edges can only be adjacent to edges labeled with~$\infty$ or~$2$. In particular, these groups include those constructed from an even FC-type Artin group by taking direct and free products with odd dihedral Artin groups.
\end{remark}

\begin{remark}
Any retraction used so far, sends a fixed generator to either a generator or to identity. We call such retractions \emph{ordinary retractions}.
\end{remark}

\begin{example}
    Not every retraction is ordinary. Indeed, consider the group:
    \[
    A_{I_2(4)} = \langle a, b\,|\, abab = baba\rangle.
    \]
    The morphism induced by~$\rho_{\{b\}} \colon A_{I_2(4)} \longrightarrow \langle b \rangle$ which maps~$a, b$ to~$b^2, b$ respectively, is a retraction but it is not an ordinary one.
\end{example}

Note that retractions used previously in \cite{antolin2015tits}, \cite{AntolinFoniqi}, and \cite{antolin2023subgroups} are ordinary. 

As a summary of this section, we have that FC type Artin groups admitting retractions, do in fact admit ordinary retractions as well. Since the presence of retractions is the main important thing on deducing results, we will focus only on the ordinary retractions.

The main results in this paper are about the class of FC type Artin groups, together with the ordinary retractions defined in \autoref{thm: characterization of FC type admitting retractions}.
However, most of the results in \Cref{sec: Retractions to standard parabolic subgroups uniquely extending to parabolic subgroups} and \Cref{sec: retractions} hold for other Artin groups admitting ordinary retractions. 

Here we provide a class of Artin groups, that together with its ordinary retractions, contains the class of FC type Artin groups with the ordinary retractions. 

\begin{definition}
We say that an Artin group~$A_\Gamma$ is \emph{(odd, odd)-free} if~$\Gamma$ does not contain two consecutive odd-labelled edges. 
\end{definition}

\begin{remark}
    Notice that our class of FC-type Artin groups with ordinary retractions is (odd, odd)-free.
\end{remark}

If we follow \autoref{thm: characterization of FC type admitting retractions}, we notice that we can try to define ordinary retractions in the same way in  (odd, odd)-free Artin groups as well. Moreover, here we obtain more groups admitting ordinary retractions than in the FC-type case.
\begin{example} Let~$\Gamma$ be the following labeled graph, with~$m, n \in \bN$.
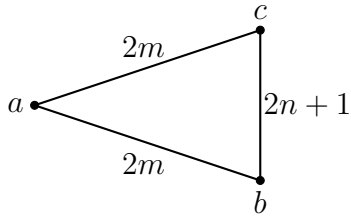
\begin{figure}[H]
\centering
\begin{tikzpicture}
%% vertices & nodes
\draw[fill=black] (0,0) circle (1.5pt) node[left] {$a$};
\draw[fill=black] (3,-1) circle (1.5pt) node[below] {$b$};
\draw[fill=black] (3,1) circle (1.5pt) node[above] {$c$};
\draw[fill=black] (2.9,0) circle (0pt) node[right] {$2n+1$};
\draw[fill=black] (1.45,0.5) circle (0pt) node[above] {$2m$};
\draw[fill=black] (1.45,-0.5) circle (0pt) node[below] {$2m$};
%edges
\draw[thick] (0,0) -- (3,1);
\draw[thick] (3,-1) -- (0,0);
\draw[thick] (3,1) -- (3,-1);
\end{tikzpicture}\label{fig: odd-odd free with retractions}
\end{figure}
The corresponding Artin group~$A_{\Gamma}$ admits ordinary retractions. Note that the~$A_{\Gamma}$ is not of FC-type if~$m > 1$.
\end{example}

\section{Retractions to standard parabolic subgroups extend uniquely to all parabolic subgroups}\label{sec: Retractions to standard parabolic subgroups uniquely extending to parabolic subgroups}

For even Artin groups, \cite{AntolinFoniqi} observed that if we have a standard retraction~$\rho_X \colon A_S \longrightarrow A_X$, then for the parabolic subgroup~$P = \alpha A_X \alpha^{-1}$, the map~$\rho_X^\alpha \colon A_S \longrightarrow P$, defined by~$\rho_X^\alpha(g) \coloneqq \alpha \rho_X(\alpha^{-1} g \alpha) \alpha^{-1}$, is also a retraction. This observation holds for any Artin group.
Indeed, for~$\gamma \in A_S$ denote by~$c_g$ the conjugation map 
$$c_{\gamma} \colon A_S \longrightarrow A_S \text{ defined by } c_{\gamma}(g) = \gamma g \gamma^{-1}.$$
Then one has the following commutative diagram:
\[
\begin{tikzcd}
A_S \arrow{r}{\rho_X^\alpha} \arrow[swap]{d}{c_{\alpha^{-1}}} & \alpha A_X \alpha^{-1}  \\
A_S \arrow{r}[swap]{\rho_X} & A_X \arrow[swap]{u} {c_{\alpha}}
\end{tikzcd}
\]
showing that~$\rho_X^\alpha$ is indeed a morphism which when restricted to~$\alpha A_X \alpha^{-1}$ gives the identity, hence it is a retraction as well.
%admitting ordinary retractions. 

Note that for the same parabolic subgroup one could have several retractions to it
% coming from retractions to different standard parabolic subgroups
, in the sense that if~$P = \alpha A_X \alpha^{-1} = \beta A_Y \beta^{-1}$, then both~$\rho_X^\alpha$ and~$\rho_Y^\beta$ are retractions to~$P$ and they could be different. In this section we will prove that this is not possible, i.e. all retractions to parabolic subgroups are determined by any associated retraction to a standard parabolic subgroup.

\begin{remark}\cite[Corollary~3.5.]{AntolinFoniqi}
In even Artin groups, a parabolic subgroup~$P=\alpha A_X \alpha^{-1}$ cannot be expressed as~$\beta A_{Y}\beta^{-1}$ with~$X\neq Y$. 
\end{remark}

\begin{lemma}
If~$P = \alpha A_X \alpha^{-1} = \beta A_{X}\beta^{-1}$ in an even Artin group~$A_S$, then~$\rho_X^\alpha = \rho_X^\beta$.
\end{lemma}
\begin{proof}
The equality~$P = \alpha A_X \alpha^{-1} = \beta A_{X}\beta^{-1}$, implies that~$A_X = \gamma A_X \gamma^{-1}$ for~$\gamma = \alpha^{-1} \beta$. In this case to show that~$\rho_X^\alpha = \rho_X^\beta$ pick~$g \in A_S$. One has the equations:
\[
\rho_X^\alpha(g) = \alpha \rho_X(\alpha^{-1} g \alpha) \alpha^{-1}, \text{ and } \rho_X^\beta(g) = \beta \rho_X(\beta^{-1} g \beta) \beta^{-1}.
\]
Rewriting~$\alpha^{-1} g \alpha$ as~$\gamma \beta^{-1} g \beta \gamma^{-1}$, and setting~$h = \beta^{-1} g \beta$ it remains to show that:
$\rho_X(\gamma h \gamma^{-1})$ is equal to~$\gamma \rho_X(h) \gamma^{-1}$.
But both of these belong to~$A_X$ as~$A_X = \gamma A_X \gamma^{-1}$, so we can apply~$\rho_X$, which yields the equality~$\gamma \rho_X(h) \gamma^{-1} = \rho_X(\gamma \rho_X(h) \gamma^{-1}) = \rho_X(\gamma h \gamma^{-1})$, as required.
\end{proof}

\smallskip
One aim in this section is to show that in FC-type Artin groups admitting retractions, the definition of~$\rho_X^\alpha$ does not depend on the presentation of~$P = \alpha A_X \alpha^{-1}$, i.e. it does not depend on the choice of~$X$ and~$\alpha$.
\begin{remark}\label{rem: conjugation on D_3}
    Differently to even Artin groups, note that in FC-type Artin groups admitting retractions one can have~$\alpha A_X \alpha^{-1} = \beta A_{Y}\beta^{-1}$ with~$X\neq Y$. Indeed, consider the group:
    \[
    A_{I_2(3)} = \langle a, b \,|\,aba = bab\rangle.
    \]
    The relation~$aba = bab$ gives the equality~$aba^{-1} = b^{-1}ab$, implying that~$aA_{\{b\}}a^{-1} = b^{-1}A_{\{b\}}b$. 
    %In other words, we have the equality~$\alpha A_X \alpha^{-1} = \beta A_{Y}\beta^{-1}$ with~$X\neq Y$, for~$\alpha = a$,~$\beta = b^{-1}$,~$X = \{b\}$, and~$Y = \{a\}$.
\end{remark}

Another goal is to find conditions for~$\rho_X^\alpha = \rho_X^\beta$ in the general setting of Artin groups. The conditions that we will obtain are sufficient but not necessary, and they hold in particular for the class of FC-type Artin groups admitting retractions.
\smallskip

First, to describe the general situation of the \autoref{rem: conjugation on D_3}, we use a lemma which is a particular case of a result of \cite{Paris} --- see also \cite[Section~2]{Cumplido2022} for an alternative explanation ---. Suppose that for the Coxeter graph~$\Gamma$ of~$A_\Gamma$ there is no edge between~$s, t$ if~$m_{s,t}=2$; then we can consider the connected components of~$\Gamma$, which we denote by~$\Delta_1,\Delta_2,\dots, \Delta_r$. In this framing, we have that~$$A_\Gamma= A_{\Delta_1}\times A_{\Delta_2}\times \cdots \times A_{\Delta_r},$$ and we call each~$A_{\Delta_i}$ an \emph{irreducible component} of~$A_{\Gamma}$. Note that the irreducible components of~$A_{\Gamma}$ are standard parabolic subgroups.

\begin{lemma}[{Adapted from \citealp[Section~4]{Paris}}]\label{lemma Luis}
Let~$A_\Gamma$ be an Artin group such that any of its irreducible components~$A_Z$ with~$|Z| > 2$ is not of spherical type, and let~$P = f A_X f^{-1}$ be a parabolic subgroup. 
%Then~$P = g A_Y g^{-1}$ with~$Y \neq X$ only if~$\langle x \rangle$ is a direct connected component of~$A_X$ and~$x$,~$x\in X$, is connected in~$\Gamma$ to some~$y$ through a path of odd labels. In this case,~$A_Y = \langle y \rangle \times A_{Y'}$ and~$A_X = \langle x \rangle \times A_{X'}$.
Then~$P = g A_Y g^{-1}$ with~$Y \neq X$ only if there exists~$x\in X$ such that~$A_X$ can be decomposed as a direct product with one factor being~$\langle x \rangle$, and~$x$ is connected in~$\Gamma$ to some~$y \in Y$ through a path of odd labels. In this case,~$A_Y = \langle y \rangle \times A_{Y'}$ and~$A_X = \langle x \rangle \times A_{X'}$, with~$X' = X \setminus \{x\}$, and~$Y' = Y \setminus \{y\}$.

\end{lemma}

Roughly speaking, an irreducible component of an Artin group can be conjugated to another irreducible component only if it has spherical type and sits inside a standard parabolic subgroup of spherical type. Also, the generator of the centre of this bigger standard parabolic subgroup needs to be the square of its Garside element. In our case, the only irreducible components that satisfy these conditions are the ones with just one generator connected through odd labels to another generator in the Coxeter graph. 

%In the even case, a parabolic subgroup~$P=gA_Xg^{-1}$ cannot be expressed as~$hA_{Y}h^{-1}$ with~$X\neq Y$. In the context of FC-type Artin groups admitting retractions, by the previous lemma, if~$X$ does not have any cyclic direct components~$\langle x\rangle$, where~$x\in X$, we are in the same situation.
%In that event, we say that~$P$ is {\it parabolic over}~$X$. Otherwise, for every cyclic direct component~$\langle x\rangle~$  of~$A_X$, since we cannot have two consecutive odd-labelled edges, there can be only one~$y$,~$x\neq y$, such that~$f\langle x \rangle f^{-1}=h\langle y\rangle h^{-1}$, for some~$f,h\in A$. This means that we have at most two choices of standard parabolic subgroup for every cyclic component. We still can say that~$P$ is parabolic over~$X$, but the choice of~$X$ is not unique. We want to prove that we can define a retraction to~$P$ and it does not depend on this choice.
\smallskip

As another auxiliary result, we describe the normalizers of standard parabolic subgroups in Artin groups such that every irreducible~$A_Z$ with~$|Z| > 2$ is not of spherical type. 

Consider~$X$ such that~$\langle x \rangle$, with~$x \in X$, is an irreducible component of~$A_X$. For~$y \in S \setminus X$, we will say that 
$$\underbrace{yxy\cdots}_{m_{x,y}-1}$$ 
and its inverse are \emph{elementary ribbons} if the irreducible component of~$A_{X \cup \{y\}}$ containing~$y$ also contains~$x$ and defines a dihedral Artin group with~$m_{x,y} > 2$. If we write~$X = \{ x \} \cup X'$, then the elementary ribbons above conjugate~$X$ to~$Y \in \{X, X' \cup \{y\}\}$; We refer to ribbon as above by the name \emph{elementary~$X$-ribbon-$Y$}.
In a more general setting, an~$X$-ribbon-$Y$ is an element~$r = r_1 \cdots r_n$, where each~$r_i$ is an elementary~$X_{i-1}$-ribbon-$X_i$, with~$X_0 = X$ and~$X_n = Y$. We denote the set of all~$X$-ribbons-$Y$ by~$\mathrm{Ribb}(X,Y)$.

\begin{conjecture}[{Adapted from \citealp[Conjecture~1]{Godelle2007}}]\label{conj Godelle}
Let~$A$ be an Artin group such that every irreducible~$A_Z$ with~$|Z| > 2$ is not of spherical type. If~$A_X$ does not have any cyclic irreducible component~$\langle x \rangle$ with~$x \in X$, then 
$$N(A_X) = A_X \times A_{X^\bot},$$ 
where 
$$X^\bot = \{y \in S \setminus X \mid xy = yx, \forall x \in X\}.$$ 
Otherwise, the elements that conjugate~$X$ to~$Y$ are given by the set 
$$\mathrm{Conj}(X,Y) = (\mathrm{Ribb}(X,Y) \cdot A_X) \cdot A_{X^\bot}.$$

\end{conjecture}

\begin{remark}\label{rem: ribbons commute}
    %Notice that, for FC-type Artin groups admitting retractions,  an odd-labelled edge can only be adjacent to an edge labelled with a~$2$ or a~$\infty$. Therefore, if~$a$ and~$b$ are connected through an edge labelled with an odd~$m$ then the only generator different from~$a$ which is a conjugate of~$a$ is~$b$ and:~$$Ribb(\{a\},\{b\})= \langle \underbrace{ba\cdots b}_{m-1}\rangle,\quad Ribb(\{a\},\{a\})= \langle \underbrace{ba\cdots a}_{m-1} \underbrace{ab\cdots b}_{m-1}\rangle.$$
Notice that, by definition, every~$X$-ribbon-$Y$ commutes with every non-cyclic irreducible component of~$X$.
\end{remark}

For every Artin group, there is a general version of ribbon and a general conjecture to describe normalizers of standard parabolic subgroup --- see \cite[Conjecture~1]{Godelle2007} --- but we will not need it for retractions in the FC case. In general, we know that FC-type Artin groups \cite[Theorem~0.3]{godelle2003parabolic} and 2-dimensional Artin groups satisfy this general ribbon conjecture \cite[Theorem~3]{Godelle2007}. Notice that our condition includes 2-dimensional Artin groups.

\begin{lemma}\label{lem triangle odd}
Let~$A$ be an Artin group admitting ordinary retractions. For every standard parabolic subgroup~$A_{\{a,b,c\}}$ with~$m_{a,b}$ and~$m_{b,c}$ odd, we have that~$m_{a,c}$ is also odd. 
\end{lemma}

\begin{proof}
Since~$\rho_{\{a,c\}}(b)\in \{a,c\}$. We have that the relation~$aca\cdots =cac \cdots$ of length~$m_{a,b}$ or~$m_{b,c}$ needs to be fulfilled in~$A_{\{a,c\}}$, which is impossible if~$m_{a,c}$ is even. 
\end{proof}

\begin{theorem}\label{thm:retractions extension}
Let~\( A \) be an Artin group admitting ordinary retractions, such that~\( A \) satisfies \autoref{conj Godelle} and every irreducible standard parabolic subgroup~\( A_Z \) with~\(|Z| > 2\) is not of spherical type. Let~\( P = fA_X f^{-1} \) be a parabolic subgroup, where~\( f \in A \). Then, the morphism
\[
\rho_X^f \colon A \longrightarrow P, \qquad \rho_X^f(a) \coloneqq f \rho_X(f^{-1} a f) f^{-1},
\]
defines a retraction. Moreover, if~\( P = gA_Y g^{-1} \) as well, then~\(\rho_X^f(a) = \rho_Y^g(a)\) for all~\( a \in A \). We denote this unique retraction by~\(\rho_P\).

\end{theorem}

\begin{proof}
 If~$fA_X f^{-1}$ is such that~$A_X$ has no cyclic irreducible component, then by \autoref{lemma Luis} we know that~$A_X$ cannot be a conjugate of a different standard parabolic subgroup. Suppose~$P = fA_X f^{-1}$ and~$P = gA_X g^{-1}$. In this case,~$g \in fN(A_X)$; hence, by \autoref{conj Godelle}, one has~$g = f \cdot \sigma \cdot c$ for some~$\sigma \in A_X$ and some~$c \in A_{X^{\bot}}$. By our definition, 
$$\rho_X^f(a) = f\rho_X(f^{-1}af)f^{-1} = f\rho_X(f)^{-1}\rho_X(a)\rho_X(f)f^{-1},$$
$$\rho_X^g(a) = g\rho_X(g^{-1}ag)g^{-1} = g\rho_X(g)^{-1}\rho_X(a)\rho_X(g)g^{-1}.$$
We have that~$g\rho_X(g)^{-1} = f \sigma c \rho_X(c)^{-1} \rho_X(\sigma)^{-1} \rho_X(f)^{-1}$. As~$\sigma \in A_X$, we have~$\rho_X(\sigma) = \sigma$; moreover,~$c$~commutes with anything in~$A_X$, so~$\rho_X(c) = 1$,~$g\rho_X(g)^{-1} = f \rho_X(f)^{-1} c$, and~$\rho_X^f(a) = \rho_{Y}^g(a)$.

\smallskip

Let~$f \langle a \rangle f^{-1}$ be a cyclic direct component of~$A_X$. 
% Since we cannot have two consecutive odd-labelled edges, there can be only one~$b$,~$a \neq b$, such that~$f A_X f^{-1} = h A_Y h^{-1}$, where~$Y = \{b\} \cup (X \setminus \{a\})$. 
Assume~$f A_X f^{-1} = h A_Y h^{-1}$. Denote the components of~$A_X$ by~$A_{X_i}$ and the components of~$A_Z$ by~$A_{Z_i}$. We need to prove that~$\rho_X^f(g) = \rho_Y^h(g)$. We can suppose that every~$A_{X_i}$ is a conjugate of~$A_{Y_i}$ by~$fh^{-1}$.

\smallskip

By \autoref{lemma Luis}, every cyclic~$A_{\{a\}}$ is conjugate by~$fh^{-1}$ to a cyclic~$A_{\{b\}}$, where~$b$ is either~$a$ or a different generator in~$Y$. For the remaining components, there is~$Y_i$ such that~$X_i = Y_i$. 
Since any~$x \in S$ is mapped to a unique component through~$\rho_X$ or~$\rho_Y$, notice that 
\[
\rho_X(f^{-1}gf) = \rho_{X_1}(f^{-1}gf) \cdots \rho_{X_r}(f^{-1}gf)
\]
and 
\[
\rho_Y(h^{-1}gh) = \rho_{Y_1}(h^{-1}gh) \cdots \rho_{Y_r}(h^{-1}gh),
\]
so the problem reduces to proving that~$h^{-1}f$ conjugates~$\rho_{X_i}(f^{-1}gf)$ to~$\rho_{Y_i}(h^{-1}gh)$.

\smallskip

Assume that~$X_i = Y_i$ and~$|X_i| > 1$. By \autoref{conj Godelle},~$h^{-1}f$~can be decomposed as the product~$f_3 f_2 f_1$, where~$f_1$ are elements commuting with~$X$,~$f_2$~is an element inside~$A_X$, and~$f_3$ are~$X$-ribbons-$Z$. Notice that~$\rho_{X_i}(f_1) = 1$,~$\rho_{X_i}(f_2)$ is the part of the direct product decomposition of~$f_2$ that lies in~$A_{X_i}$, and~$\rho_{X_i}(f_3) = 1$ (\autoref{rem: ribbons commute}). This implies
\begin{align*}
    h^{-1}f \rho_{X_i}(f^{-1}gf) f^{-1}h  & = \rho_{X_i}(h^{-1}f) \rho_{X_i}(f^{-1}gf) \rho_{X_i}(f^{-1}h) \\
    & = \rho_{X_i}(h^{-1}ff^{-1}gff^{-1}h) 
    \\
    & = \rho_{Y_i}(h^{-1}gh)
\end{align*}
as we wanted. If~$X_i = Y_i = \{a\}$, we have
$$
h^{-1}f \rho_{\{a\}}(f^{-1}gf) f^{-1}h = h^{-1}f(a^{-m'} a^m a^{m'}) f^{-1}h = a^m
$$
and
$
\rho_{\{a\}}(h^{-1}gh) = a^{-m''} a^m a^{m''} = a^m.
$
Now, suppose~$X_i = \{a\}$ and~$Y_i = \{b\}$ are different and conjugate to each other by~$h^{-1}f$. Then we have
$$
h^{-1}f \rho_{\{a\}}(f^{-1}gf) f^{-1}h = h^{-1}f(a^{-m'} a^m a^{m'}) f^{-1}h = b^m.
$$
Finally, by \autoref{lem triangle odd},~$a$~and~$b$ are connected to the same number of vertices by odd labels, so
$
\rho_{\{b\}}(h^{-1}gh) = b^{-m''} b^m b^{m''} = b^m.
$
\end{proof}

\section{Retractions and intersection of parabolic subgroups}\label{sec: Retractions and intersection of parabolic subgroups}

\label{sec: retractions}
Retractions have proven to be useful in extracting properties of parabolic subgroups of Artin groups. For instance, the following lemma holds for any Artin group, as shown in \citep*[Proposition~2.6]{moller2022parabolic}. The proof provided below, which utilizes retractions as in even Artin groups \cite[Lemma~3.4]{AntolinFoniqi}, is much simpler.

\begin{proposition}\label{lem: containing implies equality}
Let~$A_S$ be an Artin group admitting retractions,~$g,h\in A_S$ and~$X\subseteq S$.
If~$gA_Xg^{-1}\leqslant hA_Xh^{-1}$ then~$gA_Xg^{-1}=  hA_Xh^{-1}$.
\end{proposition}

\begin{proof}
One has~$gA_Xg^{-1}\leqslant hA_Xh^{-1}$ if and only if~$h^{-1}gA_Xg^{-1}h\leqslant A_X$. Now one can apply the retraction~$\rho_X$ to obtain~$h^{-1}gA_Xg^{-1}h = \rho_X(h^{-1}gA_Xg^{-1}h) = A_X$, as we wanted.
\end{proof}
\noindent In even Artin groups, retractions commute, moreover they satisfy
$$\rho_{X \cap Y}(x) = \rho_X \circ \rho_Y (x) = \rho_Y \circ \rho_X(x)$$
for every generator~$x$ (hence for every group element as well). Unfortunately, in more general cases, and specifically in the FC-case, this property does not hold. Indeed, consider the example below.
\begin{example}\label{ex: non-commuting retractions in D_3}
In the group~$A_{I_2(3)} = \langle a, b\,|\,aba = bab\rangle$, for~$X = \{a\}$, and~$Y = \{b\}$ using our ordinary retractions, we have 
$$\rho_X \circ \rho_Y (a) = a \text{ but } \rho_Y \circ \rho_X (a) = b.$$
\end{example}

However, we can identify the situations where obstructions occur. We say that an Artin group~$A_\Gamma$ is \emph{(odd, odd)-free} if~$\Gamma$ does not contain two consecutive odd-labelled edges. Notice that the condition in the proposition below is more restrictive than the one in \autoref{thm:retractions extension}, but the FC-case still fits within it.

\begin{proposition}\label{lem:retractions_modified}
Let~$A_S$ be an (odd, odd)-free Artin group admitting ordinary retractions. Let~$X, Y \subseteq S$ and~$\rho_X, \rho_Y$ the corresponding retractions. Then, for every~$s\in S$, one of the following holds:
\begin{itemize}
    \item[(i)]~$\rho_{X \cap Y}(s)=1$,
    \item[(ii)]~$\rho_{X \cap Y}(s)=\rho_X\circ \rho_Y (s)= \rho_Y\circ \rho_X(s)$.
\end{itemize}
In particular, the case (ii) always holds except when~$s \in X \setminus Y$ (or analogously~$s \in Y \setminus X$) and there is an odd-labelled edge in the Coxeter graph connecting~$s$ to~$Y \setminus X$ (respectively connecting~$s$ to~$X \setminus Y$).

\end{proposition}

\begin{proof}
We simply check all the cases coming from the definition of our retractions. We denote by~$e$ (if it exists) the odd-labelled edge with one endpoint~$s$; denote its other endpoint by~$s'$. 

\begin{itemize}
    \item \textbf{Case 1}: If~$s \in X \cap Y$, then 
  ~$s = \rho_{X \cap Y}(s) = \rho_X(s) = \rho_Y(s)$, which gives as well:
  ~$$s = \rho_{X \cap Y}(s) = \rho_X \circ \rho_Y (s) = \rho_Y \circ \rho_X(s).$$
    \item \textbf{Case 2}: Assume~$s \notin X \cup Y$. 
    \begin{itemize}
        \item If~$e$ does not exist or~$e$ is not adjacent to~$X \cup Y$, then~$\rho_X(s) = 1 = \rho_Y(s)$, which give:
      ~$$1 = \rho_{X \cap Y}(s) = \rho_X \circ \rho_Y (s) = \rho_Y \circ \rho_X(s).$$
        \item  If~$e$ is adjacent to~$X \cap Y$, then~$s' = \rho_{X \cap Y}(s) = \rho_X(s) = \rho_Y(s)$, yielding:
      ~$$s' = \rho_{X \cap Y}(s) = \rho_X \circ \rho_Y (s) = \rho_Y \circ \rho_X(s).$$
        \item  If~$e$ is adjacent to~$X \setminus Y$ (the case for~$Y \setminus X$ is analogous), then~$\rho_{X}(s) = s'$. Since we are in an (odd, odd)-free situation, we have~$\rho_{X \cap Y}(s) = 1$,~$\rho_Y(s) = 1$, and~$\rho_Y(\rho_X(s)) = \rho_Y(s') = 1$. Finally we have:
      ~$$1 = \rho_{X \cap Y}(s) = \rho_X \circ \rho_Y (s) = \rho_Y \circ \rho_X(s).$$
    \end{itemize}    
    \item \textbf{Case 3}: Suppose that~$s \in X \setminus Y$ (the case for~$s \in Y \setminus X$ is analogous), so~$\rho_X(s) = s$.
    \begin{itemize}
        \item If~$e$ is adjacent to~$X \cap Y$, then~$\rho_{X \cap Y}(s) = s'$, and~$\rho_{Y}(s) = s'$. Note that~$s' \in X \cap Y$, so~$\rho_X(s') = s' = \rho_Y(s')$. Ultimately:
      ~$$s' = \rho_{X \cap Y}(s) = \rho_X \circ \rho_Y (s) = \rho_Y \circ \rho_X(s).$$
        \item If~$e$ does not exist or is not adjacent to~$Y$, then~$\rho_Y(s) = 1$ and~$\rho_{X \cap Y}(s) = 1$
        %and~$\rho_X(s) = s$, 
        so:
      ~$$1 = \rho_{X \cap Y}(s) = \rho_X \circ \rho_Y (s) = \rho_Y \circ \rho_X(s).$$ 
        \item Finally, assume that~$e$ is adjacent to~$Y \setminus X$. This is the only case where case~$(ii)$ does not hold, as we have~$\rho_X(\rho_Y(s)) = \rho_X(s') = s$ and~$\rho_Y(\rho_X(s)) = \rho_Y(s) = s'$. \\        
        Nevertheless, we have the equality~$\rho_{X \cap Y}(s) = 1$, so case~$(i)$ holds.
    \end{itemize}
\end{itemize}
This finishes the proof.
\end{proof}

\noindent
The generators for which commutation does not occur in the previous proposition will be of special importance in this section and will form distinguished sets.

\begin{definition}\label{Definition OC}
Let~$A_S$ be an (odd, odd)-free Artin group admitting ordinary retractions. For~$X,Y\subseteq S$, denote by~$O_{X,Y}$ the elements in~$X\setminus Y$ connected to~$Y\setminus X$ through an odd labeled edge. 
%and by~$O_{Y,X}$ the elements in~$Y\setminus X$ connected to~$X\setminus Y$ through an odd label. 
We also define~$C_{X,Y}:=(X\cap Y)\cup O_{X,Y}$.
%and~$C_{Y,X}:=(Y\cap X)\cup O_{Y,X}$.
\end{definition}

%\begin{lem}[Godelle]\label{lem: proper_inclusions_parabolics}
%Let~$A$ be an FC-type Artin group generated by~$S$,~$X, Y \subseteq S$ and~$g, h \in A$. If~$g A_X g^{-1}\subsetneq A_Y$, then~$X\subsetneq Y$ and~$g\in A_Y$. In particular,~$gA_Xg^{-1} \subsetneq hA_Yh^{-1}$ implies~$X\subsetneq Y$.
%\end{lem}

\begin{remark}
    In the class of even Artin groups (consequently in right-angled Artin groups as well) we do not have odd labels, so in those groups one has~$O_{X,Y} = \emptyset$ and~$C_{X,Y} = X\cap Y$.
\end{remark}

\begin{remark}\label{rem: bijection}
    Notice that, with the (odd, odd)-free condition,~$\rho_X$ (when restricted to generators)
    %(or~$\rho_Y$) 
    defines a bijection between~$O_{X,Y}$ and~$O_{Y,X}$, and consequently between~$C_{X,Y}$ and~$C_{Y,X}$.
\end{remark}

\begin{lemma}[{\citealp{van1983homotopy}}]\label{lem: %
intersection_of_standard_parabolic_subgroups}
Let~$A_S$ be any Artin group and let~$X, Y \subseteq S$. The following equality holds:
$$A_X \cap A_Y = A_{X \cap Y}.$$
\end{lemma}

\noindent
In the following lemma, we prove that that the intersection of two parabolic subgroups over~$X$ and~$Y$ reduces to intersections of parabolic subgroups over the special sets~$C_{X,Y}$ and~$C_{Y,X}$.

\begin{remark}
    For the proof of the next lemma we want to use a convention on parabolic subgroups. For an Artin group~$A_S$ and~$X, Y \subseteq S$, we want to have~$\rho_X(A_Y) = A_{\rho_X(Y)}$. This is true under the convention that~$A_{\{1\}} \coloneqq \{1\}$ in~$A_S$, and~$A_{Z \cup \{1\}} \coloneqq A_{Z}$ in~$A_S$ for~$Z \subseteq S$. From this point onwards, we will not include the identity~$1$ in the generating set of parabolic subgroups.
\end{remark}

\begin{lemma}\label{lem: intersection_of_two_parabolic_subgroups}
Let~$A_S$ be an (odd, odd)-free Artin group admitting ordinary retractions. For any~$f, g \in A_S$ and~$X, Y \subseteq S$, there exist~$x \in A_X$ and~$y \in A_Y$ such that
$$fA_Xf^{-1} \cap gA_Yg^{-1} = fxA_{C_{X,Y}}x^{-1}f^{-1} \cap gyA_{C_{Y,X}}y^{-1}g^{-1}.$$
\end{lemma}

\begin{proof}
This proof is similar to the one of \cite[Lemma~3.3]{AntolinFoniqi}, with minor modifications to include the sets~$C_{X,Y}$ and~$C_{Y,X}$. We write it below for completeness.

\smallskip

We start with the equality
$$fA_Xf^{-1} \cap gA_Yg^{-1} = f[A_X \cap (f^{-1}g)A_Y(f^{-1}g)^{-1}]f^{-1}.$$
Set~$h = f^{-1}g$ and consider~$P = A_X \cap hA_Yh^{-1}$. Using~$P \subseteq A_X$ and~$A_X \cap A_Y = A_{X \cap Y}$ (\autoref{lem: intersection_of_standard_parabolic_subgroups}), we obtain:
\begin{align*}
P = \rho_X(P) &= \rho_X(A_X \cap hA_Yh^{-1}) \\
&\subseteq \rho_X(A_X) \cap \rho_X(hA_Yh^{-1}) \\
&= A_X \cap \rho_X(h)\rho_X(A_Y)\rho_X(h^{-1}) \\
&= \rho_X(h)A_{\rho_X(Y)}\rho_X(h)^{-1}.
\end{align*}
Notice that~$\rho_X(Y)$ consists of~$X \cap Y$ and the elements of~$X \setminus Y$ connected to~$Y \setminus X$ through an odd label. In other words,~$\rho_X(Y) = C_{X,Y}$. Setting~$x = \rho_X(h) \in A_X$, we can write the inclusion above as~$P \subseteq xA_{C_{X,Y}}x^{-1}$, and we notice that~$xA_{C_{X,Y}}x^{-1} \subseteq A_X$. Recalling as well~$P = A_X \cap hA_Yh^{-1}$, we have:
\begin{align*}
P &= (A_X \cap hA_Yh^{-1}) \cap xA_{C_{X,Y}}x^{-1} \\
&= hA_Yh^{-1} \cap (A_X \cap xA_{C_{X,Y}}x^{-1}) \\
&= hA_Yh^{-1} \cap xA_{C_{X,Y}}x^{-1}.
\end{align*}
Conjugating the last equation by~$h^{-1}$ and setting~$P' = h^{-1}Ph$ and~$k = h^{-1}x$, we obtain:
$$P' = A_Y \cap kA_{C_{X,Y}}k^{-1}.$$
Applying the same procedure as for~$P$ above, we obtain:
\begin{align*}
P' \subseteq \rho_Y(k)A_{C_{Y,X}}\rho_Y(k)^{-1},
\end{align*}
where~$C_{Y,X} = \rho_Y(C_{X,Y})$ (\autoref{rem: bijection}). Setting~$y = \rho_Y(k) \in A_Y$, we express the inclusion above as~$P' \subseteq yA_{C_{Y,X}}y^{-1} \subseteq A_Y$. Combining~$P' = A_Y \cap kA_{C_{X,Y}}k^{-1}$ and~$P' \subseteq yA_{C_{Y,X}}y^{-1}$, we have:
$$P' = (A_Y \cap kA_{C_{X,Y}}k^{-1}) \cap yA_{C_{Y,X}}y^{-1} = h^{-1}xA_{C_{X,Y}}x^{-1}h \cap (A_Y \cap yA_{C_{Y,X}}y^{-1}).$$
Using~$A_Y \cap yA_{C_{Y,X}}y^{-1} = yA_{C_{Y,X}}y^{-1}$ and~$P' = h^{-1}Ph$, we ultimately have:
$$P = xA_{C_{X,Y}}x^{-1} \cap hyA_{C_{Y,X}}y^{-1}h^{-1}.$$
Returning to our initial setup, we have~$fA_Xf^{-1} \cap gA_Yg^{-1} = fPf^{-1}$, and since~$h = f^{-1}g$, we obtain:
$$fA_Xf^{-1} \cap gA_Yg^{-1} = fxA_{C_{X,Y}}x^{-1}f^{-1} \cap gyA_{C_{Y,X}}y^{-1}g^{-1},$$
with~$x \in A_X$ and~$y \in A_Y$, as desired.
\end{proof}

%\begin{lem}
%Let~$A,B\subseteq V$ and~$f,g\in G_\Gamma$.
%Let~$H=fG_Af^{-1}$ and~$K=gG_Bg^{-1}$ be parabolic subgroups of~$G_\Gamma$. If~$H=K$ then all elements in~$A-B$ and~$B-A$ are exotic.
%\end{lem}

%\begin{proof}
%  ~$H = K$ is equivalent to~$G_A = kG_Bk^{-1}$ with~$k = f^{-1}g$. Applying~$\rho_A$ to it, we obtain~$G_A = G_{C_A}$, i.e.~$A = C_A$. Similarly,~$B = C_{B}$.
%\end{proof}

%\begin{ex}
 %   In~$G = \Gpres{a,b}{aba=bab}$ one has~$bG_{\{a\}}b^{-1} = aG_{\{b\}}a^{-1}$.
%\end{ex}

\begin{lemma}\label{lem: link exterior}
Let~$A_S$ an (odd, odd)-free Artin group admitting ordinary retractions, and let~$g \in A_S$ be any element. For~$X,Y\subseteq S$ set~$D=C_{X,Y}\cup C_{Y,X}$.
Suppose~$A_{C_{X,Y}} \cup gA_{C_{Y,X}}g^{-1}$
is not contained in a proper parabolic
subgroup of~$A_S$ and for some~$s \in S \setminus D$ one has that~$D \not\subset \lk(s)$.
Then~$A_{C_{X,Y}} \cap  gA_{C_{Y,X}}g^{-1}$
is contained in a parabolic subgroup over a
proper subset of~$C_{X,Y}$ or a proper subset of~$C_{Y,X}$.
\end{lemma}

\begin{proof}
This proof is adapted from the proof of \cite[Lemma~3.6]{AntolinFoniqi}. During this proof, we use the no-$\infty$-edge convention for the Coxeter graph. 
%the convention that the Coxeter graph associated to~$A_S$ has no edge if the associated label is~$\infty$ and a labelled edge otherwise.
Let~$P = A_{C_{X,Y}} \cap gA_{C_{Y,X}}g^{-1}$. If~$A_{S \setminus \{s\}} = gA_{S \setminus \{s\}}$, then~$g \in A_{S \setminus \{s\}}$. Consequently, both~$A_{C_{X,Y}}$ and~$gA_{C_{Y,X}}g^{-1}$ would be parabolic subgroups of~$A_{S \setminus \{s\}}$, and thus their union,~$A_{C_{X,Y}} \cup gA_{C_{Y,X}}g^{-1}$, would be contained within the proper parabolic subgroup~$A_{S \setminus \{s\}}$ of~$A_S$. This contradicts the assumptions of the proposition. 
Therefore, we assume~$A_{S \setminus \{s\}} \ne gA_{S \setminus \{s\}}$. Given this, the presentation of~$A_S$ allows us to view it as an amalgamated free product:
%\begin{equation}\label{star_splitting}
$$A_S = A_{\st(s)}*_{A_{\lk(s)}} A_{S \setminus \{s\}}.$$
%\end{equation}
Consider the Bass-Serre tree~$T$ associated with the given splitting, as discussed, for example, in \citep{DicksDunwoody}. In this tree, vertices fall into two categories: those representing left cosets of~$A_{\st(s)}$ and those representing left cosets of~$A_{S \setminus {s}}$ within~$A_S$. Vertices of different categories are the only ones that can be adjacent in~$T$.
The group~$A_S$ acts naturally on~$T$, without edge inversions. Additionally, the stabilizers of vertices in~$T$ correspond to conjugates of~$A_{\st(s)}$ or~$A_{S \setminus {s}}$, depending on the type of the vertex. The stabilizers of edges, on the other hand, are conjugates of~$A_{\lk(s)}$.
\smallskip

In the tree~$T$, the vertices~$A_{S \setminus \{s\}}$ and~$gA_{S \setminus \{s\}}$ are distinct, and yet belong to the same category. Their stabilizers are~$A_{S \setminus \{s\}}$ and~$gA_{S \setminus \{s\}}g^{-1}$, respectively. Since~$T$ is a tree, there exists a unique geodesic path~$p$ connecting~$A_{S \setminus \{s\}}$ and~$gA_{S \setminus \{s\}}$.
Since~\(s \not\in C_{X,Y} \cup C_{Y,X} \), we have:
\[
A_{C_{X,Y}} \subseteq A_{S \setminus \{s\}} \quad \text{and} \quad gA_{C_{Y,X}}g^{-1} \subseteq gA_{S \setminus \{s\}}g^{-1},
\]
which implies that the parabolic subgroups~\( A_{C_{X,Y}} \) and~\( gA_{C_{Y,X}}g^{-1} \) stabilize the vertices corresponding to the cosets~\( A_{S \setminus \{s\}} \) and~\( gA_{S \setminus \{s\}} \), respectively. The intersection~\( A_{C_{X,Y}} \cap gA_{C_{Y,X}}g^{-1} \) stabilizes the geodesic~\( p \) connecting these vertices and, consequently, any edge belonging to~\( p \).
Since the stabilizers of edges in~\( T \) are conjugates of~\( A_{\lk(s)} \), it follows that:
\[
P = A_{C_{X,Y}} \cap gA_{C_{Y,X}}g^{-1} \subseteq hA_{\lk(s)}h^{-1}
\]
for some~\( h \in A_S \).
We can express~\( P \) as:
\begin{align}\label{eq: big intersection of P}
\nonumber
    P & = A_{C_{X,Y}} \cap gA_{C_{Y,X}}g^{-1} \cap hA_{\lk(s)}h^{-1} \\ 
    & = (A_{C_{X,Y}} \cap hA_{\lk(s)}h^{-1}) \cap (gA_{C_{Y,X}}g^{-1} \cap hA_{\lk(s)}h^{-1}).
\end{align}

We follow \autoref{Definition OC} and consider the sets~\( C_{D,\lk(s)} \) and~\( C_{\lk(s),D} \). \\
{\bf Claim.} There exists some~\( y \in D \) such that~\( y \not\in C_{D,\lk(s)} \). \\
{\it Proof of Claim.}
Assume, for the sake of contradiction, that every~\( y \in D \setminus \lk(s) \) is odd-connected to some~\( y' \in \lk(s) \setminus D \). Now focus on the Artin group~\( A_{\{y,y',s\}} \). Since the edge~\( \{y, y'\} \) has an odd label, then by the (odd, odd)-free condition, the edge~\( \{s, y'\} \) has an even label (note that~$y' \in \lk(s)$, so this edge cannot be~\( \infty \)). Thus, by \autoref{lem:impossible triangles infinity}, the edge~\( \{s, y\} \) cannot be labeled with~\( \infty \), as otherwise the Artin group~\( A_{\{y,y',s\}} \) would not admit (ordinary) retractions; implying that~$y \in \lk(s)$. This implies that every~\( y \in D \) would need to be in~\( \lk(s) \), contradicting the assumption that~\( D \subsetneq \lk(s) \).
\smallskip

By \autoref{lem: intersection_of_two_parabolic_subgroups}, the intersection in \autoref{eq: big intersection of P}  can be expressed as an intersection of parabolic subgroups over the sets
$$C_{C_{X,Y}, \lk(s)},\; C_{\lk(s), C_{X,Y}},\; C_{C_{Y,X}, \lk(s)}, \text{ and }~C_{\lk(s), C_{Y,X}},$$
all of which are subsets of~\( D \).

%By the claim from the previous paragraph and \autoref{lem: intersection_of_two_parabolic_subgroups}, the intersection of a conjugate of~\( A_{D} \) with a conjugate of~\( A_{\lk(s)} \) is contained within a proper parabolic subgroup of~\( D \), specifically, a conjugate of~\( A_{D \setminus \{y\}} \).

By the claim above there is some~\( y \in D \) such that~\( y \not\in C_{D,\lk(s)} \). Recall from our notation that~$D=C_{X,Y}\cup C_{Y,X}$, and so
$C_{D,\lk(s)} = C_{C_{X,Y},\lk(s)} \cup C_{C_{Y,X},\lk(s)}$.
Therefore, depending on where the~\( y \) of the claim lies
%(recall that~$y \in D$ and~$y \notin C_{D,\lk(s)}$),
either~\( A_{C_{X,Y}} \cap hA_{\lk(s)}h^{-1} \) is contained in a parabolic subgroup over a proper subset of~\( C_{X,Y} \), or~\( gA_{C_{Y,X}}g^{-1} \cap hA_{\lk(s)}h^{-1} \) is contained in a parabolic subgroup over a proper subset of~\( C_{Y,X} \), as desired.
\end{proof}

%\begin{lem}\label{lem: reducing subgraphs} \textcolor{blue}{This lemma is weaker than Lemma~\ref{lem: proper_inclusions_parabolics}. We don't need it.}
%Let~$\Delta$ be a subgraph of~$\Gamma$,~$A\subseteq V_\Delta$ and~$g,t\in G_\Gamma$. 
%If~$G_A\cup gG_Ag^{-1}$ is contained in~$tG_\Delta t^{-1}$,
%then there is~$h\in G_\Delta$ such that:
%\begin{enumerate}
%\item[(i)]~$G_A=hG_Ah^{-1}$ if and only if~$G_A=gG_Ag^{-1}$. 
%\item[(ii)]~$G_A\cap hG_Ah^{-1}$ is~$\Delta$-parabolic if and only if~$G_A\cap gG_Ag^{-1}$ is~$\Gamma$-parabolic. 

%\item[(iii)]~$G_A\cap hG_Ah^{-1}$ is contained in a~$\Delta$-parabolic over a proper subset of~$A$  if and only if~$G_A\cap gG_Ag^{-1}$ is~$\Gamma$-parabolic over a proper subset of~$A$.
%\end{enumerate} 
%\end{lem}

%\textcolor{red}{Note from Islam: one needs to show that the following conditions is satisfied for our groups}
\subsection{Artin groups with property~$\cC$}\label{subsec: Artin groups with property C}
\noindent
We say that a family of (odd, odd)-free Artin groups has \emph{property~$\cC$} if it satisfies the following conditions:
\begin{enumerate}
    \item The family consists of Artin groups whose underlying Artin graphs admit retractions and is closed under taking subgraphs.
    \item For every Artin group~$A_S$ in the family, for any subsets~$X, Y \subseteq S$, and any element~$g \in A_\Gamma$, if all vertices~$x \in S \setminus (C_{X,Y} \cup C_{Y,X})$ satisfy~$\st(x) = S$, then either:
    \begin{enumerate}
        \item~$A_{C_{X,Y}} = g A_{C_{Y,X}} g^{-1}$, or
        \item~$A_{C_{X,Y}} \cap g A_{C_{Y,X}} g^{-1} \leqslant d A_H d^{-1}$ for some subset~$H$ properly contained in~$C_{X,Y}$ or~$C_{Y,X}$, and some~$d \in A_S$.
    \end{enumerate}
\end{enumerate}

\begin{remark}
    \cite[Theorem 5.2]{AntolinFoniqi} shows that even Artin groups of FC type satisfy the condition~$\mathcal{C}$.
\end{remark}

\begin{lemma}\label{lemma restandardisation}
Let~$A_S$ be an (odd, odd)-free Artin group that admits ordinary retractions,
and let~$X, Y \subseteq S$. If~$A_{C_{X,Y}} \cup gA_{C_{Y,X}}g^{-1}$ is contained in a proper parabolic subgroup~$h A_Z h^{-1}$ of~$A_S$, then there exist subsets~$X', Y' \subseteq Z$ and elements~$f_1, f_2 \in A_Z$ such that
$$
h^{-1}(A_{C_{X,Y}} \cap gA_{C_{Y,X}}g^{-1})h = f_1 A_{C_{X',Y'}} f_1^{-1} \cap f_2 A_{C_{Y',X'}} f_2^{-1}.
$$
\end{lemma}

\begin{proof}
We will use the result of \cite[Theorem~1.1]{BlufsteinParis} - which states that if~$uA_Ru^{-1} \leqslant A_Q$ for some~$u \in A_S$ and~$R, Q \subseteq S$, then~$uA_Ru^{-1} = vA_{Q'}v^{-1}$ for some~$v \in A_Q$, and some~$Q' \subseteq Q$.
Now we know that~\( h^{-1} A_{C_{X,Y}} h = t_1 A_{X'} t_1^{-1} \) for some~\( X' \subseteq Z \) and~\( t_1 \in A_Z \). Similarly,~\( h^{-1} g A_{C_{Y,X}} g^{-1} h = t_2 A_{Y'} t_2^{-1} \) for some~\( Y' \subseteq Z \) and~\( t_2 \in A_Z \). Now a direct application of \autoref{lem: intersection_of_two_parabolic_subgroups} for the intersection~$t_1 A_{X'} t_1^{-1} \cap t_2 A_{Y'} t_2^{-1}$ gives the appropriate elements~$f_1, f_2$, providing the desired result.
\end{proof}

The following theorem reduces the problem of intersection of parabolic subgroups to property~$\cC$.

\begin{theorem}\label{thm: reduction to stars}
Let~$\mathcal{F}$ be a family of (odd, odd)-free Artin groups. If~$\mathcal{F}$ satisfies condition~$\mathcal{C}$, and~$A_S$ belongs to~$\mathcal{F}$, then the intersection of two parabolic subgroups of~$A_S$ is again a parabolic subgroup. 
\end{theorem}

\begin{comment}
\begin{theorem}\label{thm: reduction to stars}
If a family of (odd, odd)-free Artin groups~$A_S$ admitting ordinary retractions has property~$\mathcal{C}$, then the intersection of two parabolic subgroups of~$A_S$ is again a parabolic subgroup. 
\end{theorem}
\end{comment}

\begin{proof}
This is an adaptation of \cite[Proposition~3.8]{AntolinFoniqi}. During this proof, we consider the no-$\infty$-edge convention for Coxeter graphs. 
%we consider the Coxeter graph~$\Gamma$ associated to~$A_S$ with no edges if the associated label is~$\infty$ and a labelled edge otherwise. 
Let~\( A_\Gamma = A_S \) be an Artin group as described in the statement, and let~\( P, Q \) be two parabolic subgroups of~\( A_\Gamma \). By \autoref{lem: intersection_of_two_parabolic_subgroups}, we can assume there exist~\( X, Y \subseteq S \) and~\( h_1, h_2 \in A_\Gamma \) such that~\( P \cap Q = h_1 A_{C_{X, Y}} h_1^{-1} \cap h_2 A_{C_{Y, X}} h_2^{-1} \). Therefore,~\( P \cap Q \) is a parabolic subgroup if and only if~\( A_{C_{X, Y}} \cap gA_{C_{Y, X}}g^{-1} \) is a parabolic subgroup, where~\( g = h_1^{-1}h_2 \).

\smallskip

If~\( A_{C_{X, Y}} \cup gA_{C_{Y, X}}g^{-1} \) is contained in a proper parabolic subgroup of~\( A_S \), then by \autoref{lemma restandardisation} we can replace~\( \Gamma \) with a proper subgraph~\( \Delta \) and replace~\( g \) with some~\( h \in A_\Delta \). Note that~\( A_{\Delta} \) also belongs to~$\mathcal{F}$.
%remains in the class of (odd, odd)-free Artin groups that admit ordinary retractions. 
Therefore, we can revert to the initial notation and further assume that~\( A_{C_{X, Y}} \cup gA_{C_{Y, X}}g^{-1} \) is not contained in a proper parabolic subgroup of~\( A_\Gamma \).

\smallskip

We aim to show that for any~\( X, Y \subseteq S \) and any~\( g \in A_\Gamma \), the intersection~\( A_{C_{X, Y}} \cap gA_{C_{Y, X}}g^{-1} \) is a parabolic subgroup. Let~\( D \coloneqq C_{X, Y} \cup C_{Y, X} \). Our proof is by induction on~\( |D| \).
If~\( |D| = 0 \), then~\( |C_{X, Y}| = |C_{Y, X}| = 0 \), i.e.,~\( A_{C_{X, Y}} \) and~\( A_{C_{Y, X}} \) are trivial, and the result follows. We now assume that~\( |D| > 0 \) and that the result holds for parabolic subgroups over smaller sets.

\smallskip

If there exists~\( x \in S \setminus D \) such that~\( D \) is not contained in~\(\lk(x)\), then by \autoref{lem: link exterior}, we have~\( A_{C_{X,Y}} \cap gA_{C_{Y,X}}g^{-1} \leqslant dA_Hd^{-1} \) for some~\( H \) properly contained in~\( C_{X,Y} \) or~\( C_{Y,X} \), and some~\( d \in A_\Gamma \). Suppose, without loss of generality, that~\( H \subsetneq C_{X,Y} \) --- the other case works similarly up to conjugacy ---. Now, additionally, we can assume that~$d \in A_{C_{X,Y}}$; indeed, we obtain it, by applying~$\rho_{C_{X, y}}$ to~\( A_{C_{X,Y}} \cap gA_{C_{Y,X}}g^{-1} \leqslant dA_Hd^{-1} \).
Hence,
\begin{align*}
    A_{C_{X,Y}} \cap gA_{C_{Y,X}}g^{-1} & = A_{C_{X,Y}} \cap gA_{C_{Y,X}}g^{-1} \cap d A_H d^{-1} \\
    & = d A_H d^{-1} \cap gA_{C_{Y,X}}g^{-1}.
\end{align*}

%By the inductive hypothesis and \cite[Theorem~1.1]{BlufsteinParis}, we know that~\( A_{C_{X,Y}} \cap d A_H d^{-1} = d'A_{H'}d'^{-1} \) for some~\( H' \subsetneq C_{X,Y} \) and~\( d' \in A_{C_{X,Y}} \). 
\noindent
Since~\( |H| \cup C_{Y,X}| < |D| \), the desired result follows by induction.
\smallskip

So, from now on, we assume that for all~\( x \in S \setminus D \), one has~\( D \subseteq \lk(x) \). We now argue by induction on~\( N = \sharp \{ x \in S \setminus D : \st(x) \neq S \} \).
In the case~\( N = 0 \), since we have property~\(\cC\), we have that either~\( A_{C_{X,Y}} = gA_{C_{Y,X}}g^{-1} \) and hence~\( A_{C_{X,Y}} \cap gA_{C_{Y,X}}g^{-1} \) is parabolic, or~\( A_{C_{X,Y}} \cap gA_{C_{Y,X}}g^{-1} \leqslant dA_H d^{-1} \) for some~\( H \subsetneq D \). In the latter case, we can use the same argument as before to prove that~\( A_{C_{X,Y}} \cap gA_{C_{Y,X}}g^{-1} \) is a parabolic subgroup.
Now, assume that~\( N > 0 \) and the result is known for smaller values of~\( N \).

\smallskip
Since~\( N > 0 \) and~\( D \subseteq \lk(x) \) for all~\( x \in S \setminus D \), there exist~\( a, b \in S \setminus D \) that are not linked by an edge. Setting~\( Z = \st(a) \),~\( U = S \setminus \{a\} \), and~\( L = \lk(a) \), we obtain an amalgamated free product:
\[
A = A_Z *_{A_L} A_U,
\]
and there is an associated Bass-Serre tree~\(\mathcal{T}\) corresponding to this splitting. Consider the edges~\( A_L \) and~\( gA_L \) on~\(\mathcal{T}\). Let~\( A_L, g_1A_L, \ldots, g_nA_L, gA_L \) be a sequence of edges in the unique geodesic in~\(\mathcal{T}\) connecting~\( A_L \) and~\( gA_L \).
If~\( A_L = gA_L \) (i.e.,~\( n = 0 \)) then~$g \in A_L$; taking into account that~\( D \subseteq L \), we have that~\( A_{C_{X,Y}} \cup g A_{C_{Y,X}}g^{-1} \) is contained in~\( A_L \), which is a proper parabolic subgroup of~\( A_S \), which contradicts our hypothesis.

\smallskip

We will proceed by induction on~\( n \). Fix an~$n$ and suppose that the result is true for all smaller values. By the construction of~\(\mathcal{T}\), for any~\( i = 0, \ldots, n \), one has either~\( g_i^{-1}g_{i+1} \in A_Z \) or~\( g_i^{-1}g_{i+1} \in A_U \), where~\( g_0 = 1 \) and~\( g_{n+1} = g \).
The intersection~\( A_{C_{X,Y}} \cap gA_{C_{Y,X}} g^{-1} \) stabilizes the endpoints of the geodesic path in~\(\mathcal{T}\), hence it stabilizes the entire path. As the stabilizer of a geodesic in a tree is the intersection of the stabilizers of its edges, we have the equality
\[
A_{C_{X,Y}} \cap gA_{C_{Y,X}} g^{-1} = A_{C_{X,Y}} \cap g_1A_L g_1^{-1} \cap \ldots \cap g_nA_L g_n^{-1} \cap gA_{C_{Y,X}} g^{-1}.
\]

By \autoref{lem: intersection_of_two_parabolic_subgroups}, applied to~${A_{C_{X,Y}}}\cap g_i A_Lg_i^{-1}$, and the fact that~${C_{C_{X,Y},L}}\subseteq L$ --- recall that~$D\subset L$ --- we have that there are~$z_i\in A_L$ such that
$$A_{{C_{X,Y}}}\cap g_i A_{L} g_i^{-1} \text{ is equal to } A_{{C_{X,Y}}}\cap g_iz_i A_{{C_{X,Y}}} z_i^{-1}g_i^{-1}.$$ 
Note that~$(g_iz_i)^{-1}(g_{i+1}z_{i+1}) = z_i^{-1}(g_i^{-1}g_{i+1})z_{i+1}$. By a slight abuse of notation we can replace~$g_iz_i$ by~$g_i$ and still have that~$g_i^{-1}g_{i+1} \in A_{Z}$ or~$g_i^{-1}g_{i+1} \in A_{U}$, for any~$i = 0, \ldots, n$ where~$g_0 = 1$ and~$g_{n+1} = g$ (because~$z_i, z_{i+1}$ belong in~$A_L$). 
Hence:
$$A_{{C_{X,Y}}}\cap gA_{{C_{Y,X}}} g^{-1} = A_{{C_{X,Y}}}\cap g_1A_{{C_{X,Y}}}g_1^{-1} \cap \ldots \cap g_nA_{{C_{X,Y}}}g_n^{-1}\cap gA_{{C_{Y,X}}} g^{-1}.$$
Consider now~$$g_iA_{{C_{X,Y}}}g_i^{-1} \cap g_{i+1}A_{{C_{X,Y}}}g_{i+1}^{-1} = g_i( A_{{C_{X,Y}}} \cap  g_{i}'A_{{C_{X,Y}}}g_{i}'^{-1})g_i^{-1},$$ with~$g'_i=g_i^{-1}g_{i+1}$. As~$g_i' \in A_{Z}$ or~$g_i' \in A_{U}$, the intersection~$A_{C_{X,Y}}\cap g_i'A_{C_{X,Y}}g_i'^{-1}$ takes place either in~$A_Z$ or in~$A_U$. 
Note that the number of vertices in~$Z \setminus D$ (respectively in~$U\setminus D$) whose star is not equal to~$Z$ (respectively~$U$) is smaller than~$N$. Therefore, by the induction hypothesis on~$N$, the intersection is either equal to~$A_{C_{X,Y}}$ or is contained in a proper parabolic subgroup over~$C_{X,Y}$. 

Suppose first that~$C_{X,Y}=C_{Y,X}$. If for all~$i$ we get~$A_{C_{X,Y}}\cap g_i'A_{C_{X,Y}}g_i'^{-1} = A_{C_{X,Y}}$ then~$A_{C_{X,Y}}\cap gA_{C_{Y,X}} g^{-1} = A_{C_{X,Y}}$, otherwise,~$A_{C_{X,Y}}\cap gA_{C_{Y,X}} g^{-1}$ is contained in a proper parabolic subgroup of~$A_{C_{X,Y}}$, and we have the result by induction on~$|D|$. 

Otherwise, for~$C_{X,Y} \neq C_{Y,X}$ one has~$|C_{X,Y}|<|D|$ and we also have the same conclusion by induction on~$|D|$.
\end{proof}

If the class of FC-type Artin groups admitting retractions would satisfy condition~$\mathcal{C}$ then the intersection of any two parabolic subgroups, would be again parabolic. However, this article has not managed to answer whether condition~$\mathcal{C}$ is satisfied by FC-type Artin groups admitting retractions. Nevertheless, the work presented here provides some new results --- see the discussion on the following remark. 

\begin{remark}
In the general class of FC-type Artin groups we know that the intersection of a spherical-type parabolic subgroup and any other parabolic subgroup is again a spherical parabolic subgroup (\citealp[Theorem~3.1]{Rose} and \citealp[Corollary~1.2]{moller2022parabolic}). We do not know if this is true for any two parabolic subgroups of non-spherical type. Thanks to the results in this paper, even if we do not know if Condition~$\cC$ is satisfied, we can provide some partial answers to this problem. For instance, we know now that if~$X$ and~$Y$ have trivial intersection and no~$x\in X$ is connected to~$Y$ through an odd-label edge, then~$A_X\cap g^{-1} A_Y g$ is the trivial subgroup for any~$g$ in the group. Also, we know that if~$X$ and~$Y$ intersect in a irreducible component~$A_Z$ of the group and~$O_{X,Y}=O_{Y,X}=\emptyset$, then~$A_X\cap g^{-1} A_Y g = A_Z$.

For example, consider the FC-type Artin group~$A$ corresponding to the following Coxeter graph with the no-$2$-edge convention:% where the 2-labels have been omitted:

\bigskip

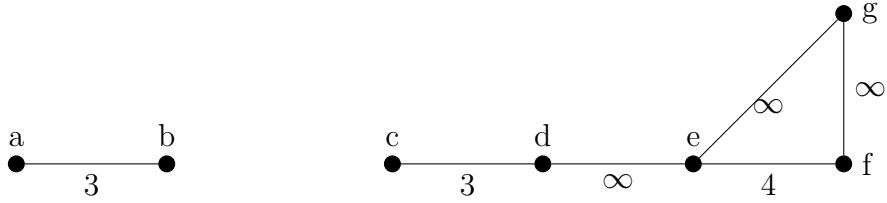
\begin{center}
\begin{tikzpicture}

\tikzstyle{vertex}=[circle, draw, fill=black, inner sep=0pt, minimum size=6pt]

\node[vertex] (A1) at (0,0) {};
\node[vertex] (B1) at (2,0) {};
\draw (A1) -- node[below] {3} (B1);

\node[above=0.1cm] at (A1) {a};
\node[above=0.1cm] at (B1) {b};

\node[vertex] (A2) at (5,0) {};
\node[vertex] (B2) at (7,0) {};
\node[vertex] (C2) at (9,0) {};
\node[vertex] (D2) at (11,0) {};
\node[vertex] (E2) at (11,2) {};

\draw (A2) -- node[below] {3} (B2);
\draw (B2) -- node[below] {$\infty$} (C2);
\draw (C2) -- node[below] {4} (D2);
\draw (C2) -- node[below] {$\infty$} (E2);
\draw (D2) -- node[right] {$\infty$} (E2);

\node[above=0.1cm] at (A2) {c};
\node[above=0.1cm] at (B2) {d};
\node[above=0.1cm] at (C2) {e};
\node[right=0.1cm] at (D2) {f};
\node[right=0.1cm] at (E2) {g};

\end{tikzpicture}
\end{center}

\bigskip

We discuss the intersection of pairs of some non-spherical parabolic subgroups~$A_X$, and~$g^{-1}A_Y g$; denote this intersection as~$H=A_X\cap g^{-1}A_Y g$.  If~$X=\{d,e\}$ and~$Y=\{f,g\}$, we know that~$H=\{1\}$. The same happens for the following pairs of generating sets~$(X,Y)$: 
\begin{align*}
  \begin{array}{lr}
 (\{a,b,c,d,e\}, &\{f,g\}) \\
 (\{a,b,d,e\}, & \{f,g\})    \\
 (\{a,c,d,e\}, &\{f,g\}) \\
 (\{a,d,e\}, & \{f,g\})    \\
 (\{b,c,d,e\}, &\{f,g\})\\
 \end{array} & &
 \begin{array}{lr}
 (\{b,d,e\}, & \{f,g\})    \\
 (\{c,d,e\}, & \{f,g\}) \\
 (\{c,d,e\}, &\{a,f,g\}) \\
 (\{c,d,e\}, &\{a,b,f,g\}) \\
 (\{c,d,e\}, &\{b,f,g\})
\end{array}
\end{align*}

%~$$\begin{array}{lr}
%  (\{a,b,c,d,e\}, &\{f,g\}) \\
%  (\{a,b,d,e\}, & \{f,g\})    \\
%  (\{a,c,d,e\}, &\{f,g\}) \\
%  (\{a,d,e\}, & \{f,g\})    \\
%  (\{b,c,d,e\}, &\{f,g\})\\
%  (\{b,d,e\}, & \{f,g\})    \\
%  (\{c,d,e\}, & \{f,g\}) \\
%  (\{c,d,e\}, &\{a,f,g\}) \\
%  (\{c,d,e\}, &\{a,b,f,g\}) \\
%  (\{c,d,e\}, &\{b,f,g\})
% \end{array}$$
 If~$X=\{a,b,d,e\}$ and~$Y=\{a,b,g,f\}$, or if~$X=\{a,b,c,d,e\}$ and~$Y=\{a,b,g,f\}$, then~$H=A_{a,b}$. Moreover, consider the following pairs~$(X,Y)$:

\begin{align*}
  \begin{array}{lr}
 (\{a,b,c,d,e\}, &\{a,f,g\}) \\
 (\{a,b,c,d,e\}, &\{b,f,g\}) \\
 (\{a,b,d,e\}, & \{a,f,g\})    \\
 (\{a,b,d,e\}, & \{b,f,g\})    \\
 (\{a,c,d,e\}, &\{a,f,g\}) \\
 (\{a,c,d,e\}, &\{a,b,f,g\}) \\
 (\{a,c,d,e\}, &\{b,f,g\}) \\
 (\{a,d,e\}, &\{a,f,g\}) \\
 \end{array} & &
  \begin{array}{lr}
 (\{a,d,e\}, &\{a,b,f,g\}) \\
 (\{a,d,e\}, &\{b,f,g\}) \\
 (\{b,c,d,e\}, &\{a,f,g\})\\
 (\{b,c,d,e\}, &\{a,b,f,g\})\\
 (\{b,c,d,e\}, &\{b,f,g\})\\
 (\{b,d,e\}, &\{a,f,g\})\\
 (\{b,d,e\}, &\{a,b,f,g\})\\
 (\{b,d,e\}, &\{b,f,g\})
\end{array}
\end{align*}

% ~$$\begin{array}{lr}
%  (\{a,b,c,d,e\}, &\{a,f,g\}) \\
%  (\{a,b,c,d,e\}, &\{b,f,g\}) \\
%  (\{a,b,d,e\}, & \{a,f,g\})    \\
%  (\{a,b,d,e\}, & \{b,f,g\})    \\
%  (\{a,c,d,e\}, &\{a,f,g\}) \\
%  (\{a,c,d,e\}, &\{a,b,f,g\}) \\
%  (\{a,c,d,e\}, &\{b,f,g\}) \\
%  (\{a,d,e\}, &\{a,f,g\}) \\
%  (\{a,d,e\}, &\{a,b,f,g\}) \\
%  (\{a,d,e\}, &\{b,f,g\}) \\
%  (\{b,c,d,e\}, &\{a,f,g\})\\
%  (\{b,c,d,e\}, &\{a,b,f,g\})\\
%  (\{b,c,d,e\}, &\{b,f,g\})\\
%  (\{b,d,e\}, &\{a,f,g\})\\
%  (\{b,d,e\}, &\{a,b,f,g\})\\
%  (\{b,d,e\}, &\{b,f,g\})
% \end{array}$$
In these cases we know that~$H$ is the intersection of two parabolic subgroups over subsets of~$\{a,b\}$. Since~$A_{a,b}$ is an irreducible component of the group, the problem 
reduces to an intersection of groups in~$A_{\{a,b\}}$, so we know that~$H$ is a parabolic subgroup by \citep[Theorem~9.5]{cumplido2019parabolic}.  
\end{remark}

\section{Coherence for the FC case}
\label{sec: Coherence}
A group is called \emph{coherent} if all its finitely generated subgroups are also finitely presented. In this work, we provide an explicit characterization of the coherence of our groups, highlighting their similarities with other classes of Artin groups that admit retractions --- e.g. the class of RAAGs~\citep{droms1987graph}, which are groups only having 2 and~$\infty$ edges, and even Artin groups of FC type~\citep{antolin2023subgroups} ---.

A (simplicial) graph~$\Gamma$ is called {\it chordal} if all cycles of four or more vertices have a {\it chord}, which is an edge that is not part of the cycle but connects two vertices of the cycle.
Droms provided the following simple characterization of coherent RAAGs. We will consider the Coxeter graphs with the no-$\infty$-edge convention. %without edge when the label of the edge is~$\infty$.

\begin{theorem}[{\citealp[Theorem~1]{droms1987graph}}]
Let~$A$ be a RAAG based on~$\Gamma$.~$A$~is coherent if and only if~$\Gamma$ is chordal.
\end{theorem}

For a simplicial labeled graph~$\Gamma$, denote by~$\Gamma^{\leq 2}$ the graph obtained from~$\Gamma$ by removing all the edges with label~$> 2$. In the same spirit as in \citep[Theorem~6.2]{antolin2023subgroups} we have the following:
\begin{theorem}\label{thm: coherence}
Let~$A_\Gamma$ be an Artin group of FC type admitting retractions.
Then~$A_\Gamma$ is coherent if and only if both~$\Gamma$ and~$\Gamma^{\leq 2}$ are chordal.
\end{theorem}

This theorem can be deduced from the classification of coherent Artin groups, which was settled in \citep{Wise2013} after the reduction in \citep{gordon2004artin}, see \citep[Section~3]{zearra2024around} for a summary, and the lemma below.

\begin{figure}[H]
\centering
\begin{tikzpicture}

%% vertices & nodes
\draw[fill=black] (0,0) circle (1.5pt) node[left] {$a$};
\draw[fill=black] (2,-1) circle (1.5pt) node[below] {$b$};
\draw[fill=black] (4,0) circle (1.5pt) node[right] {$c$};
\draw[fill=black] (2,1) circle (1.5pt) node[above] {$d$};
\draw[fill=black] (1.9,0) circle (0pt) node[right] {$m$};

\draw[fill=black] (6,0) circle (0pt) node[right] {with~$m > 2$};
%%% edges
\draw[thick] (0,0) -- (2,1);
\draw[thick] (2,1) -- (4,0);
\draw[thick] (4,0) -- (2,-1);
\draw[thick] (2,-1) -- (0,0);
\draw[thick] (2,1) -- (2,-1);
\end{tikzpicture}\caption{The corresponding Artin group is not coherent as~$\langle a^2, b^2, c^2, d^2 \rangle \simeq F_2 \times F_2$.}\label{fig: non-coherent 4-vertex situation}
\end{figure}
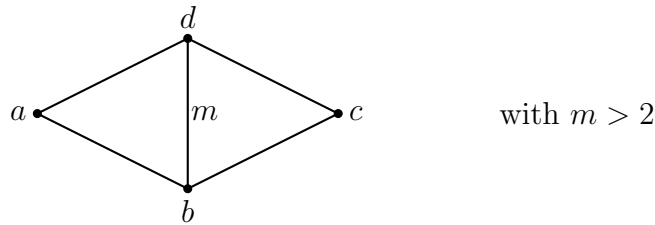

\begin{lemma}[{\citealp[Theorem~3.1]{zearra2024around}}]\label{lem: coherence of Artin groups}
Let~$\Gamma$ be a labelled finite simplicial graph. The Artin group~$A_{\Gamma}$ based on~$\Gamma$ is coherent if and only if the following three conditions hold:
\begin{itemize}[itemsep = 0pt]
\item[(i)]~$\Gamma$ is chordal;
\item[(ii)] for any complete subgraph~$\Delta \subseteq \Gamma$ of 3 or 4 vertices,~$\Delta$ has at most one edge with a label different to 2;
\item[(iii)]~$\Gamma$ has no subgraph as in \autoref{fig: non-coherent 4-vertex situation}.
\end{itemize}
\end{lemma}

\begin{proof}[Proof of \autoref{thm: coherence}]
First suppose that conditions of \autoref{lem: coherence of Artin groups} are satisfied. We have that~$\Gamma$ is chordal, so it remains to show that~$\Gamma^{\leq 2}$ is chordal as well. Assume for contradiction that~$\Gamma^{\leq 2}$ is not chordal, and so it contains an induced cycle~$C$ of length at least~$4$. 
If the length of~$C$ is bigger than~$4$, then~$\Gamma$ being chordal implies that there are two chords in~$\Gamma$ with label bigger than~$2$ and an edge of~$C$ forming a triangle, contradicting the fact that we can have only~$(2, 2, k)$-triangles with~$k \geq 2$.
So, assume then that the induced cycle~$C$ has length~$4$ in~$\Gamma^{\leq 2}$. Since~$\Gamma$ cannot contain subgraphs appearing in \autoref{fig: non-coherent 4-vertex situation}, we conclude that vertices of~$C$ form a complete graph in~$\Gamma$.
This means that~$A_C$ is isomorphic to~$A_{I_2(s)}\times A_{I_2(t)}$ with~$s,t>2$, which contains~$\mathbb{F}_2\times \mathbb{F}_2$ as a subgroup, which is not coherent. 

\smallskip

For the converse, suppose that both~$\Gamma$ and~$\Gamma^{\leq 2}$ are chordal. In particular condition~$(i)$ in \autoref{lem: coherence of Artin groups} holds. Condition (iii) is also immediate as otherwise~$\Gamma^{\leq 2}$ would not be chordal.
Now for condition (ii) assume first that~$\Delta$ is a triangle; in this case obviously we can have at most one edge with label greater than~$2$. Instead, if~$\Delta$  is complete of~$4$ vertices, with more than one edge with label greater than~$2$, then~$\Delta$ has two non-adjacent edges with labels greater than~$2$. But this means that~$\Gamma^{\leq 2}$ is not chordal.
\end{proof}

\bigskip
\noindent{\textbf{{Acknowledgments}}} 
This article is part of the research project PID2022-138719NA-I00, financed by MCIN/AEI/10.13039/501100011033/FEDER, Spain, UE. Mar\'ia Cumplido was also financed by an individual Ram\'on y Cajal 2021 Fellowship. Islam Foniqi acknowledges support from the EPSRC Fellowship grant EP/V032003/1 ‘Algorithmic, topological and geometric aspects of infinite groups, monoids and inverse semigroups’. Bruno Cisneros was also financed by CONAHCYT Research Grant Ciencia de Frontera 2019 CF 217392.
The authors acknowledge support from The London Mathematical Society (Collaborations with Developing Countries - Scheme 5, 52305).

\bibliography{main}

\begin{thebibliography}{23}
\providecommand{\natexlab}[1]{#1}
\providecommand{\url}[1]{\texttt{#1}}
\expandafter\ifx\csname urlstyle\endcsname\relax
  \providecommand{\doi}[1]{doi: #1}\else
  \providecommand{\doi}{doi: \begingroup \urlstyle{rm}\Url}\fi

\bibitem[Antol{\'\i}n and Foniqi(2023)]{antolin2023subgroups}
Y.~Antol{\'\i}n and I.~Foniqi.
\newblock Subgroups of even {A}rtin groups of {FC}-type.
\newblock arXiv preprint arXiv:2305.17292, 2023.

\bibitem[Antol{\'\i}n and Minasyan(2015)]{antolin2015tits}
Y.~Antol{\'\i}n and A.~Minasyan.
\newblock Tits alternatives for graph products.
\newblock \emph{Journal f{\"u}r die reine und angewandte Mathematik (Crelles
  Journal)}, 2015\penalty0 (704):\penalty0 55--83, 2015.

\bibitem[Antolín and Foniqi(2022)]{AntolinFoniqi}
Y.~Antolín and I.~Foniqi.
\newblock Intersection of parabolic subgroups in even {Artin} groups of
  {FC}-type.
\newblock \emph{Proceedings of the Edinburgh Mathematical Society}, 65\penalty0
  (4):\penalty0 938–957, 2022.

\bibitem[Blufstein and Paris(2023)]{BlufsteinParis}
M.~Blufstein and L.~Paris.
\newblock Parabolic subgroups inside parabolic subgroups of {A}rtin groups.
\newblock \emph{Proceedings of the American Mathematical Society}, 2023.

\bibitem[Bourbaki(1968)]{Bourbaki}
N.~Bourbaki.
\newblock \emph{\'El\'ements de math\'ematique. Fasc. XXXIV. Groupes et
  alg\`ebres de Lie. Chapitre IV: Groupes de Coxeter et syst\`emes de Tits.
  Chapitre V: Groupes engendr\'es par des r\'eflexions. Chapitre VI: syst\`emes
  de racines.}
\newblock Number 1337 in Actualit\'es Scientifiques et Industrielles. Hermann,
  Paris, 1968.

\bibitem[Brieskorn and Saito(1972)]{BrieskornSaito}
E.~Brieskorn and K.~Saito.
\newblock {Artin-Gruppen und Coxeter-Gruppen}.
\newblock \emph{Inventiones Mathematicae}, 17:\penalty0 245--271, 1972.

\bibitem[Charney and Davis(1995)]{CharneyDavis}
R.~Charney and M.~W. Davis.
\newblock The $k(\pi, 1)$-problem for hyperplane complements associated to
  infinite reflection groups.
\newblock \emph{Journal of the American Mathematical Society}, 8\penalty0
  (3):\penalty0 597--627, 1995.

\bibitem[Concini and Salvetti(2000)]{ConciniSalvetti}
C.~Concini and M.~Salvetti.
\newblock {Cohomology of Coxeter groups and Artin groups}.
\newblock \emph{Mathematical Research Letters}, 7:\penalty0 213--232, 03 2000.

\bibitem[Cumplido(2022)]{Cumplido2022}
M.~Cumplido.
\newblock The conjugacy stability problem for parabolic subgroups of {Artin}
  groups.
\newblock \emph{Mediterrean Journal of Mathematics}, 19\penalty0 (19):\penalty0
  237, 2022.

\bibitem[Cumplido et~al.(2019)Cumplido, Gebhardt, Gonz{\'a}lez-Meneses, and
  Wiest]{cumplido2019parabolic}
M.~Cumplido, V.~Gebhardt, J.~Gonz{\'a}lez-Meneses, and B.~Wiest.
\newblock On parabolic subgroups of {A}rtin--{T}its groups of spherical type.
\newblock \emph{Advances in Mathematics}, 352:\penalty0 572--610, 2019.

\bibitem[Dehornoy and Paris(1999)]{DehornoyParis}
P.~Dehornoy and L.~Paris.
\newblock {Gaussian Groups and Garside Groups, two Generalisations of Artin
  Groups}.
\newblock \emph{Proceedings of the London Mathematical Society}, 79\penalty0
  (3):\penalty0 569--604, 1999.

\bibitem[Dicks and Dunwoody(1989)]{DicksDunwoody}
W.~Dicks and M.~J. Dunwoody.
\newblock \emph{Groups acting on graphs}, volume~17 of \emph{Cambridge Studies
  in Advanced Mathematics}.
\newblock Cambridge University Press, Cambridge, 1989.

\bibitem[Droms(1987)]{droms1987graph}
C.~Droms.
\newblock Graph groups, coherence, and three-manifolds.
\newblock \emph{J. Algebra}, 106\penalty0 (2):\penalty0 484--489, 1987.

\bibitem[Godelle(2003)]{godelle2003parabolic}
E.~Godelle.
\newblock Parabolic subgroups of {A}rtin groups of type {FC}.
\newblock \emph{Pacific Journal of Mathematics}, 208\penalty0 (2):\penalty0
  243--254, 2003.

\bibitem[Godelle(2007)]{Godelle2007}
E.~Godelle.
\newblock {Artin–Tits groups with CAT(0) Deligne complex}.
\newblock \emph{Journal of Pure and Applied Algebra}, 208\penalty0
  (1):\penalty0 39–52, 2007.

\bibitem[Gordon(2004)]{gordon2004artin}
C.~Gordon.
\newblock Artin groups, 3-manifolds and coherence.
\newblock \emph{Boletín de la Sociedad Matemática Mexicana}, 10:\penalty0
  193--198, 2004.

\bibitem[Lopez~de Gamiz~Zearra and
  Mart{\'\i}nez~P{\'e}rez(2024)]{zearra2024around}
J.~Lopez~de Gamiz~Zearra and C.~Mart{\'\i}nez~P{\'e}rez.
\newblock Around subgroups of {A}rtin groups: derived subgroups and
  acylindrical hyperbolicity in the even {FC}-case.
\newblock arXiv preprint arXiv:2405.16641, 2024.

\bibitem[M\"{o}ller et~al.(2023)M\"{o}ller, Paris, and
  Varghese]{moller2022parabolic}
P.~M\"{o}ller, L.~Paris, and O.~Varghese.
\newblock On parabolic subgroups of {A}rtin groups.
\newblock \emph{Israel Journal of Mathematics}, 261\penalty0 (2):\penalty0
  809–840, 2023.

\bibitem[Morris-Wright(2021)]{Rose}
R.~Morris-Wright.
\newblock {Parabolic subgroups in FC-type Artin groups}.
\newblock \emph{Journal of Pure and Applied Algebra}, 225\penalty0
  (1):\penalty0 106468, 2021.

\bibitem[Paris(1997)]{Paris}
L.~Paris.
\newblock {Parabolic Subgroups of Artin Groups}.
\newblock \emph{Journal of Algebra}, 196\penalty0 (2):\penalty0 369--399, 1997.

\bibitem[Paris(2014)]{Paris2014}
L.~Paris.
\newblock {$K(\pi ,1)$ conjecture for {Artin} groups}.
\newblock \emph{Annales de la Facult\'e des sciences de Toulouse :
  Math\'ematiques}, Ser. 6, 23\penalty0 (2):\penalty0 361--415, 2014.

\bibitem[{Van der Lek}(1983)]{van1983homotopy}
H.~{Van der Lek}.
\newblock \emph{The homotopy type of complex hyperplane complements}.
\newblock PhD thesis, Katholieke Universiteit te Nijmegen, 1983.

\bibitem[Wise(2013)]{Wise2013}
D.~T. Wise.
\newblock The last incoherent {A}rtin group.
\newblock \emph{Proceedings of the American Mathematical Society}, 141\penalty0
  (1):\penalty0 139--149, 2013.

\end{thebibliography}
%\bibliographystyle{plain}
%\end{thebibliography}

\bigskip
\noindent\textit{\\ Bruno Aarón Cisneros de la Cruz,\\
CONAHCYT - Unidad Oaxaca del Instituto de Matemáticas de la UNAM\\ 
León No. 2, Col. Centro, Oaxaca de Juárez, Oaxaca (Mexico)\\}
{email: bruno@im.unam.mx}

\noindent\textit{\\ Mar\'ia Cumplido,\\
Departamento de Álgebra, Facultad de Matemáticas, Universidad de Sevilla\\ 
Avenida de la Reina Mercedes s/n, 41009, Seville (Spain)\\}
{email: cumplido@us.es}

\noindent\textit{\\ Islam Foniqi,\\
The University of East Anglia\\ 
Norwich (United Kingdom)\\}
{email: i.foniqi@uea.ac.uk}
\end{document}